\documentclass[11pt]{article}
\input{epsf.sty}
\usepackage{hyperref}
\usepackage{amsmath, amssymb}
\newtheorem {thm}{Theorem}[section]
\newtheorem {prop}[thm]{Proposition} 
\newtheorem {lem}[thm]{Lemma}
\newtheorem {cor}[thm]{Corollary}
\newtheorem {defn}[thm]{Definition}

\def\Cox{\hfill \Box}

\def\Z{{\Bbb Z}}
\def\R{{\Bbb R}}
\def\P{{\Bbb P}}

\def\E{{\Bbb E}}
\def\f{{\varphi}}
\def\0{{\bf 0}}

\def\ba{{\backslash}}
\def\sb{{\subset}}

\def\D{\Delta}
\def\a{\alpha}
\def\sm{\setminus}
\def\ba{\setminus}

\def\d{\delta}

\def\f{\varphi}
\def\phi{\varphi}

\def\g{\gamma}

\def\td{\tilde}

\def\l{\lambda}

\def\s{\sigma}

\def\z{\zeta}

\def\D{\Delta}

\def\L{\Lambda}

\def\O{\Omega}

\def\P{\Phi}

\def\T{\T}

\def\GG{{\cal G}}

\newcommand{\triplenorm}{| \hspace{-0.3mm} | \hspace{-0.3mm} |}

\begin{document}
\title{The Posterior metric and \\  the Goodness
 of Gibbsianness \\
for transforms of Gibbs measures}

%\thanks{Work partially
%supported by   MURST (2004-06),  Cofin:  Prin 2004028108.}
%\thanks{Work partially
%%supported by   ?? }}
\author{
Christof K\"ulske
%\thanks{Research supported by Deutsche Forschungsgemeinschaft}
\footnote{ University of Groningen, Institute  of Mathematics and
Computing Science, Postbus 407 , 9700  AK Groningen, The
Netherlands
%EURANDOM, LG 1.34,
\texttt{kuelske@math.rug.nl}, \texttt{
http://www.math.rug.nl/$\sim$kuelske/ }}\, and
Alex A. Opoku%\thanks{Research supported by XXXXXXX}
\footnote{University of Groningen, Department of Mathematics and
Computing Science, Postbus 407 , 9700  AK Groningen, The
Netherlands , \texttt{A.opoku@math.rug.nl }}
%\footnote{
%University of Groningen,
%????????????
%Centre for Theoretical Physics,
%Nijenborgh 4,
%9747AG Groningen
%The Netherlands
%\texttt{aenter@phys.rug.nl},
%\texttt{a.c.d. sdf},
%\texttt{http://www..... }}
}

\maketitle

\begin{abstract} We present a general method to derive
continuity estimates for conditional probabilities of general
(possibly continuous) spin models subjected to local
transformations. Such systems arise in the study of a stochastic
time-evolution of Gibbs measures or as noisy observations.

We exhibit the minimal necessary structure for such double-layer
systems. Assuming no a priori metric on the local state spaces, we
define the posterior metric on the local image space.
 We show that it allows in a natural way to divide
the local part of the continuity estimates from the spatial part
(which is treated by Dobrushin uniqueness here). We show in the
concrete example of the time evolution of 
rotators on the $(q-1)$-dimensional sphere how this
method can be used to obtain estimates in terms of the familiar
Euclidean metric.

 \end{abstract}



\vfill\eject

\section{Introduction} \label{sect:intro}

The absence or presence of phase transitions lies at the heart of
mathematical statistical mechanics of equilibrium systems. A phase
transition in an order parameter that can be directly observed is of
an obvious interest for the system under investigation. Moreover
sometimes also the presence or absence of phase transitions is
linked in a more subtle way to the properties of the system under
investigation. In fact, it is understood that "hidden
phase transitions" in an internal system that is not directly
observable are responsible for the failure of the Gibbs property for 
a variety of important measures that appear as transforms of
different sorts of Gibbs measures. For the mechanisms of how to
become non-Gibbs and background on renormalization group type of
pathologies and beyond, see the reviews \cite{ACDR, ACD, DEZ}.

Now, the first part of the analysis of an interacting system begins
with an understanding of the "weak coupling regime" and proving
results based on absence of phase transitions when the system
variables behave as a perturbation of independent ones. There is a
variety of competing ways to our disposition to do so, giving
related but usually not equivalent results, notably Dobrushin's
uniqueness theory \cite{DOB,GOR}, expansion methods, and percolation
and coupling methods.

Indeed, when it works, Dobrushin uniqueness has a lot of advantages,
being not very technical, but very general, requiring little explicit
knowledge of the system and providing explicit estimates on decay of
correlations. Moreover, it implies useful properties generalizing
those of independent variables. As an example of such a useful property
we mention Gaussian concentration estimates of functions of the
system variables which are obtained as a corollary when there is an
estimate on the Dobrushin interaction matrix available
\cite{KUL3,KUL4}. Especially when we are talking about continuous
spin systems a Dobrushin uniqueness approach seems favorable, since cluster
expansions are often applicable only with some technical effort
\cite{ZAH1,KUL6}, and percolation and coupling are not directly
available.

%When possible it is here recommendable to expand around a discrete
 %spin system in some sense.

A particular interest has been in recent times in the study of the
loss and recovery of the Gibbs property of an initial Gibbs measure
under a stochastic time-evolution. The study started in \cite{ACD1}
where the authors focussed on the  evolution of a Gibbs measure of an
Ising model under high-temperature spin-flip Glauber dynamics.
 The main phenomenon observed therein was the loss
of the Gibbs property after a certain transition time when the
system was started at an initial low temperature state. The measure
stays non-Gibbs forever when the initial external field was zero.
More complicated transition phenomena between Gibbs and non-Gibbs
are possible at intermediate times when there is no spin-flip
symmetry: The Gibbs property is recovered again at large but finite
values of time in the presence of non-vanishing external magnetic
fields in the external measure. A complete analysis of the
corresponding Ising mean-field system in zero magnetic field was
given in \cite{KUL2} where the authors analyzed the time-temperature
dynamic phase diagram describing the Gibbs non-Gibbs transitions. In
the analysis also the phenomenon of symmetry breaking in the set of
bad configurations was detected, meaning that a bad configuration
whose spatial average does not preserve the spin flip symmetry of the model
appears.

What remains of these phenomena for continuous spins? The case of
site-wise independent diffusions of continuous spins on the lattice
starting from the Gibbs-measure of  a special double-well potential
was considered in \cite{KUL5}. It was shown therein that a similar
loss of Gibbsianness will occur if the initial double-well potential
is deep enough. In contrast to the Ising model, this loss however is
a loss without recovery, so the measure stays non-Gibbs for all
sufficiently large times. This is due to the unbounded nature of the
spins. Short-time Gibbsianness is proved to hold also in this model.
While these results hold for a continuous spin model, the method of
proof is nevertheless based on the investigation of a "hidden
discrete model", exploiting the particular form of the
Gibbs-potential. In \cite{WIO} the authors  studied  models for 
compact spins, namely the planar rotor models on the circle subjected
to diffusive time-evolution. It is shown therein that starting with 
an initial low-temperature Gibbs measure, the time-evolved measure 
obtained for infinite- or high-temperature dynamics stays Gibbs for 
short times and for the corresponding  initial infinite- or high-
temperature Gibbs measure under infinite- or high-temperature 
dynamics, the time-evolved measure stays Gibbs forever. Their 
analysis uses the machinery of cluster expansions, as earlier developed in 
\cite{ROE}. Even before it was shown that the whole process of space-time 
histories can be viewed as a Gibbs measure\cite{DEU}. This is interesting in itself, 
but does not imply that fixed-time projections are Gibbs. 

Short-time Gibbsianness in all  these models follows from
uniqueness of a hidden or internal system. While this is expected to
hold very generally, results that are not restricted to particular
models appear only for discrete spin systems \cite{LEN}. The present
paper now narrows the gap. It provides a proof of the preservation
of the Gibbs property of the time-evolved Gibbs-measures of a
general
 continuous spin system under site-wise independent dynamics,
 for short times,
 even when the initial measure is in the strong coupling regime.
 More generally than for time-evolution,
 we prove our results directly for general two layer systems,
consisting of (1) a Gibbs-measure in the first layer, that is (2)
subjected to local transition kernels mapping the first layer
variables to second layer variables. This generalizes the notion of
a hidden Markov model where the second layer plays the role of a
noisy observation. Such models have motivation in a variety of
fields. Let us mention for example that they appear in biology as models of gene
regulatory networks   where the vertices of the network are genes and
the variables model gene expression activity.

A measure is a Gibbs measure when the single-site conditional
probabilities depend on the conditioning in an essentially local
way. Our main statement (Theorem \ref{main1}) is an explicit upper
 bound on
the continuity of the single-site conditional probabilities of the second
layer system as a function of the conditioning. This is valid when the transition kernels don't
fluctuate too much, even when the first layer system is in a strong
coupling regime. Our result holds for discrete or continuous compact
state spaces and general interactions and is based on Dobrushin
uniqueness. To formulate the resulting continuity estimate for the
conditional probabilities we don't need any a priori metric
structure on the local spin spaces: The natural metric on the second
layer single spin space is created by the variational distance
between the a-priori measures in the first layer that are obtained
by conditioning on second layer configurations(see Theorem
 \ref{main1}).

On the way to this result, we exhibit a simple criterion for
Dobrushin uniqueness for Gibbs-measures (of one layer).  It is easy
to check and can be of use beyond the study of (non)-Gibbsianness.

Intuitively,  it demands that the sum over the interaction terms in
the Hamiltonian coupling the sites $i$ and $j$  should not fluctuate
too much when it is viewed as a random variable at the site $i$
under the corresponding local a-priori measure (see Definition
\ref{def1}). So even when one has a large interaction, better
concentration properties of the a priori measures can still imply an
overall  small Dobrushin constant. This is a generalization of the
simple large-field criterion ensuring Dobrushin-uniqueness in the
Ising model ( see p.147 example 8.13 of \cite{GOR} and \cite{ISR}) to general
 spaces
 (Theorem \ref{Dob3}).

In Theorem \ref{concen} we state as a corollary that "concentration
 implies
concentration". By this we mean  that there are Gaussian
concentration inequalities for functions of the coupled system, with
explicit decay rate (even when there is strong coupling) if the a
priori measures concentrate well enough.

 The criterion we need for the study of the second layer model is based on  the
 description of the interplay between the possible largeness of the
initial interaction and the strength of the coupling to the second
layer found in Theorem \ref{Dob3} (when the initial apriori measures
  are
replaced with   conditional apriori measures ). To ensure Gibbsianness of the second layer model, we
thus need small fluctuations of the initial Hamiltonian w.r.t. the
a-priori measures in the first layer that are obtained by
conditioning on second layer configurations. The estimates on the
spatial memory of the single-site second layer conditional
probabilities follow naturally by evoking Dobrushin-uniqueness
estimates on comparison of the Gibbs-measures with perturbed
specifications and chain-rule type of arguments.

To illustrate the simplicity of our approach to get explicit
estimates on the spatial decay we prove short-time Gibbsianness of (Heisenberg)-model of
$(q-1)$-dimensional rotators for general $q\geq2$ under diffusive
time-evolution on the $(q-1)$-spheres, and provide an explicit estimate on
the time-interval for which the time-evolved measure stays Gibbs.
This will be supplemented by arguments that are more specific 
to the rotators which give us precise continuity estimates in terms 
of the Euclidean distances on the spheres. 

In Section 2 we formulate our main results. In Section 3 we provide
the proofs of Theorem \ref{Dob3} and \ref{concen}, in Section 4 we
provide the proof of Theorem \ref{main1}, and in Section 5 we
provide the proofs of Theorem \ref{thm:main-intro} and  Proposition
\ref{Festimates} and provide some related results. We also give the proof of Theorem \ref{fuzzy}
 in  Section 5.

\section{Main Results}

\subsection{A criterion for Dobrushin uniqueness for concentrated a
 priori measures}

Let  $G$ be a countable vertex set, and assume that $\s=(\s_i)_{i\in
G}$ are spin-variables taking values in a measurable (standard
Borel) space $S$ (single-spin space). In our general setup we don't
need to make a metric structure on $S$ explicit. We further denote by
$\O=S^G$ the configuration space of our system equiped with the Borel
$\s$-algebra.

Let $\g$ be the Gibbs specification(collection of finite-volume
conditional distributions that depend in a  continuous way on the
conditioning ) for a given interaction potential
$\Phi=(\Phi_A)_{A\sb G}$ (where $\P_A:S^G\mapsto \mathbb{R}$ are
functions that depend only on the spin-variables in $A$ for  finite
subsets $A$ of $G$) and a priori probability measure $\a$ on the
single-site spaces, i.e for any finite $\L\sb G $ and $\bar\s\in
S^G$ we define $\g_\L(\cdot|\bar\s)\in\g$ as

\begin{eqnarray}\label{specifi}
\g_\L(d\s_\L|  { \bar\s ):= \exp\bigl( -\sum_{A\cap \L\neq
\emptyset}\Phi_A(\s_\L\bar\s_{G\ba\L} }
)\bigr)\prod_{i\in\L}\a(d\s_i)/Z_\L( \bar \s)
\end{eqnarray}
with the normalization constant $Z_\L(   \bar \s)$.

We assume the summability property
\begin{eqnarray}\label{triple}
||| \P|||:=\sup_{i\in G}\sum_{A\ni i}|A|\Vert
\Phi_{A}\Vert_{\infty}<\infty
\end{eqnarray}
for the interaction $\P$. In the sequel we will always write $i$ for
$\{i\}$, $i^c$ for $G\ba \{i\}$ and $\L^c$ for $G\ba\L$. We further
denote by $C=(C_{ij})_{i,j\in G}$  the {\it Dobrushin
interdependence matrix}, with entries given by
\begin{eqnarray}\begin{split}\label{Ci}
C_{ij}=\sup_{\zeta,\eta\in
 \Omega;\;\zeta_{j^c}=\eta_{j^c}}\Vert
 \g_i(\cdot|\zeta)-\g_i(\cdot|\eta)\Vert.
\end{split}
\end{eqnarray}
where $\Vert
\nu_1-\nu_2\Vert:=\sup_{f;|f|\leq1}|\nu_1(f)-\nu_2(f)|=\frac{1}{2}
\l(|h_1-h_2|)$ whenever  $\nu_1$ and $\nu_2$ are  probability
measures that are absolutely continuous with respect to the measure
$\l$  with $\l$-densities $h_1$ and $h_2$ respectively (i.e. $\Vert
\nu_1-\nu_2\Vert$ is one half of the variational distance between
$\nu_1$ and $\nu_2$ ). The corresponding \textit{Dobrushin constant}
is also given as
\begin{equation*}
c:=\sup_{i\in G}\sum_{j\in G}C_{ij}.
\end{equation*}
and we recall that whenever $c<1$ ( Dobrushin uniqueness condition)
then $\g$ admits at most one   Gibbs measure  \cite{DOB,
GOR}. It is known that for a potential $\P$ satisfying
(\ref{triple}) there is a sufficiently small $\beta$ such that $
\beta \Phi$ satisfies Dobrushin uniqueness and the measure is in a
small coupling regime. We will prove Dobrushin uniqueness for a potential
with possibly very large (but finite) (\ref{triple}) when the measure 
$\a$ concentrates. In fact, we can also deduce Dobrushin uniqueness
for weak coupling from the bound we will provide on the Dobrushin's contant
$c$.

For our purposes we employ the following definition.

\begin{defn} \label{def1} For a function $F:S^G\mapsto \R$ we define
the {\bf $\a;i,j$-deviation} $\text{dev}_{\a;i,j}$ of $F$ to be
 \begin{equation}\begin{split}\label{deval}
\text{dev}_{\a;i,j}(F):=\sup_{{\zeta,\eta\in
S^G}\atop{\z_{j^c}=\eta_{j^c}}}\inf_{B}\int  \a(d\s_i)\Bigr|
F(\s_i\eta_{i^c})-F(\s_i\zeta_{i^c})-B\Bigl|.
\end{split}
\end{equation}
\end{defn}

This quantity is the worst-case linear deviation of the variation of
$F$ at the site $j$ viewed as a random variable w.r.t. to $\s_i$
under $\a(d\s_i)$. Note that clearly the deviation is bounded by
$\d_j(F)$ the $j$th {\it oscillation} of $F$ , i.e.\\
$\text{dev}_{\a;i,j}(F)\leq \d_j(F)= \sup_{{\zeta,\eta\in
S^G}\atop{\z_{j^c}:=\eta_{j^c}}}\Bigr| F(\eta)-F(\zeta)\Bigl|$.

Then our first result is as follows.

\begin{thm}\label{Dob3}
The Dobrushin constant $c$ is bounded by
\begin{equation}\begin{split}\label{c-bound}
c\leq \sup_{i\in G}\sum_{j\in G}\exp\Bigr(
\sum_{A\supset\{i,j\}}\delta(\Phi_A)\Bigl)\text{dev}_{\a;i,j}(H_{i})\cr
\end{split}
\end{equation}
 where $\d(\P_A)$  is the oscillation of $\P_A$ defined as
  $\d(\P_A):=\sup_{w,\bar w\in S^G}|\P_A(w)-\P_A(\bar w)|$ and
 $H_V:=\sum_{A\cap V\neq \emptyset}\Phi_A$
\end{thm}
\bigskip

The use of this criterion lies in the fact that, even when the
interaction potential is large, $\text{dev}_{\a;i,j}(H_{i})$ can be
small, when $\a$ is close to a Dirac measure. A simple example for
this to happen is an Ising model at large external field. As a less
trivial application of the criterion to a spin-model where $S$ is
not discrete we discuss the Gauss-Weierstrass kernel in the rotator
example of Section \ref{Short-Time} where we prove short-time
Gibbsianness.

Of course, when the potential is small to begin with, the r.h.s. 
of (\ref{c-bound}) will be small, independently of $\a$, so the theorem can 
be used for both strong couplings and concentrated a priori-measures 
and weak coupling.

\subsection{Concentration implies concentration}

Dobrushin uniqueness implies also the existence of a Gibbs measure
$\mu$ (if the local spin space $S$ is standard Borel (Theorem 8.7.
\cite{GOR}).) This unique measure $\mu$ then has further nice
properties; e.g. general Gaussian estimates on the concentration of
an observable $F$ around its mean hold \cite{KUL3, KUL4}.  We 
believe that the  concentration result below  is worth mentioning.

The estimate on the Dobrushin matrix $C$ that leads to the upper
bound (\ref{c-bound}) on $c$  then also implies the Gaussian
concentration estimate which we will  give in Theorem \ref{concen}.

%We further denote $C_t$ by the
%transpose of $C$  and   $c_t=sup_{j\in G}\sum_{i\in G}C_{i,j}$ by
%the Dobrushin constant of $C_t$.

\begin{defn} We call the matrix $B$ with entries
\begin{equation}\begin{split}
B_{ij}&:= \text{dev}_{\a;i,j}(H_{i}). \cr
%s&:=\exp\Bigl(\sup_{i\neq j}\sum_{A\supset\{i,j\} }\d(\Phi_A)\Bigr)
\end{split}
\end{equation}
the deviation matrix of the potential $\Phi$ w.r.t. $\alpha$.
\end{defn}
To formulate the concentration theorem let us write
$||B||_1:=\sup_{j\in G}\sum_{i\in G}B_{i,j}$ and
$||B||_\infty:=\sup_{i\in G}\sum_{j\in G}B_{i,j}$ for the
corresponding matrix-norms.

\begin{thm}\label{concen}
 Suppose $||B||_1,||B||_\infty <\frac{1}{s}$ where
$s:=\exp\Bigl(\sup_{i\neq j}\sum_{A\supset\{i,j\}
}\d(\Phi_A)\Bigr)$. Then for  any bounded measurable function
$F(\s)$
 and
 $\forall r\geq 0$ holds the inequality
\begin{align}\mu \Bigl(F\left(\s\right)-
\mu \bigl(F\left(\s\right)\bigr)\geq r \Bigr)\leq \exp\left(-
\Bigl(1-s\Vert B\Vert_\infty\Bigr)\Bigl(1-s\Vert B\Vert_1\Bigr)
\frac{r^2}{2 \bigl\Vert\underline\d(F)\bigr\Vert_{l^2}^2}
\right)\end{align} Here we have written
$\bigl\Vert\underline\d(F)\bigr\Vert_{l^2}^2\equiv \sum_{i\in
G}(\d_{i}(F))^2$.
  \end{thm}

\subsection{Two-layer models - Goodness of Gibbsianness}
\label{twolayer}

Let us now formulate our assumptions on a two-layer system over a
graph $G$. To each vertex will be associated two local state spaces.
A particular example will be given by the site-wise independent
time-evolution of Section \ref{Short-Time}. So, in general let $S$
and additionally $S'$ be measurable (standard Borel) spaces. This
implies in particular existence of all regular conditional
probabilities. Again, no a priori metric will be used explicitly. We
refer to $S$ as the initial (first layer) spin space  and to $S'$ as
the image (second layer) spin space. Let the {\it joint a priori
measure} $K(d\s_i,d\eta_i)$ be a Borel probability measure on the
product space $S\times S'$. We assume non-nullness of $K$
(positivity of measure for all open sets). We assume further that
$K$ can be written in the form $K(d\s_i,d\eta_i)=k(\s_i,\eta_i)\a(d\s_i)
\a'(d\eta_i)$ where $\a(d\sigma_{i})\equiv \int_{S'}K(d\sigma_{i},d\eta_i)$
and $\a'(d\eta_{i}) \equiv \int_{S}K(d\sigma_{i},d\eta_i)$ with $k>0$.

Our initial model (probability measure on $S^G$)  is by definition a
Gibbs distribution for the specification given in terms of the
potential $\Phi$ according to (\ref{specifi}) where we now put as an
a priori measure the marginal of $K$ on the first layer, that is
$\a(d\sigma_{i})\equiv \int_{S'}K(d\sigma_{i},d\eta_i)$ . It is
important to note that we don't assume uniqueness of the
Gibbs measure for this specification. In practice $\a$ might be
given beforehand and $K$ is then obtained by specifying a transition
kernel $K(d\eta_i|\s_i)$ from the first layer to the second layer.
We will always  denote by $\s_i\in S$ the  local variable (spin) for
the initial model and $\eta_i\in S'$ the  local variable (spin) for
the image model.

Let $\mu(d\s)$ be a Gibbs measure for the first layer for potential
$\Phi$ and a priori measure $\a$. Our aim is then: Study the
conditional probabilities of the second layer measure defined by
\begin{equation*}\label{eq:aim}
\begin{split}
&\mu'(d\eta):=\int_{S^G}\mu(d\s)\prod_{i\in G}K(d\eta_i|\s_i).
\end{split}
\end{equation*}

This form appears for instance in the study of a stochastic
time evolution,
 starting from an initial measure $\mu$ where the kernel
$K(d\eta_i|\s_i)$ will be dependent on time and is applied
independently over the spins (infinite-temperature dynamics). In
case studies it has been observed that the map $\mu\mapsto \mu'$ may
create an image measure that is not a Gibbs measure anymore. On the
other hand, in all examples observed, Gibbsianness was preserved at
short times where $K_t$  is a small perturbation of $\d_{\eta_i}(d\s_i)$. 
We aim here to give a criterion that implies this in
all generality, not using any specifics of the model but only the
relevant underlying structure. In particular we are not restricting
ourselves to discrete spin spaces.

Our main result Theorem  \ref{main1} is a criterion for  the Gibbs
property of the second layer measure that is easily formulated and
verified in concrete examples. Moreover, we give explicit bounds on
the dependence of the conditional probabilities of the second layer
measure on the variation of the conditioning.

We said that we will not use any a priori metric on the spaces $S$
and $S'$; indeed the natural metric that shall be used for
continuity in this setup shall be given by the variational distance
of the conditional a priori measures in the first layer, conditional
on the second layer.

\begin{defn}  We call
\begin{equation*}\label{naturalmetric}
\begin{split}
&d'(\eta_j,\eta'_j):=\Vert \a_{\eta_j}- \a_{\bar \eta_j}\Vert
\end{split}
\end{equation*}
the posterior (pseudo-)metric associated to $K$ on the second layer
space.

Here $\a_{\eta_i}(d\s_i)=K(d\s_i|\eta_i)$ are the a priori measures
in the first layer that are obtained by conditioning on second layer
configurations.
\end{defn}

 In the language of statistics, $\a_{\eta_i}$ is the "posterior measure"
 depending
on the observation $\eta_i$ in the second layer single spin space. 
Stated abstractly, the metric $d'$ is
the pull back-metric of the map $\eta_i\mapsto \a_{\eta_i}(d\s_i)$ 
from single-site configurations in the second layer to single-site measures
in the first layer. 
While this metric seems to be non-explicit, we will show in the
rotator example how it can be estimated in terms of a more familiar
metric (Euclidean metric).

It is well-known that an investigation of the Gibbs property of the
second layer measure must be based on an analysis of the first layer
conditional on configurations in the second layer \cite{ACD1, DEZ,
 ACD}. So, our estimates will naturally
contain quantities that reflect this aspect. The main ingredient
will be a matrix $\bar B$ that is a uniform bound (over possible
conditionings) on the conditional deviation matrix $B(\eta)$ of the
first layer system. More precisely, let us put
\begin{equation}\begin{split}
\bar C_{ij} &:=\exp\Bigr(
\sum_{A\supset\{i,j\}}\delta(\Phi_A)\Bigl)\bar B_{ij}\quad \text{
where }\cr \bar B_{ij}&:=
\sup_{\eta_i}\text{dev}_{\a_{\eta_i};i,j}(H_{i}).\cr
\end{split}
\end{equation}
We warn the reader not to confuse
$\text{dev}_{\a_{\eta_i};i,j}(H_{i})$ with
$\text{dev}_{\a;i,j}(H_{i})$. While the second quantity may be big
and correspondingly the unconstrained first layer system in a non-uniqueness
regime, the first one might still be small and correspondingly 
the constrained layer system in a uniqueness regime.  This is e.g. the case for a
time-evolution started at low temperature, for small times. We denote by $\g'$
the class of all finite-volume conditional distributions of the transfored model
with full $\eta$-conditioning. Then we
have the following theorem.

\begin{thm}\label{main1} Suppose that the first layer system has an
 infinite-volume Gibbs measure $\mu=\lim_{n}\mu^{\bar \s}_{\L_n}$
obtained for a boundary condition $\bar \s$ and along a suitable
 sequence of volumes $\L_n$.

Suppose  further that $\sup_i\sum_{j}\bar C_{ij}<1$.

\begin{enumerate}
\item
Then $\g'$ is a specification and the second layer measure $\mu'$ is a Gibbs measure for the 
specification $\g'$.
 \item  $\g'$ satisfies the continuity estimate
\begin{equation}\begin{split}\label{deno}
&\Big\Vert\g'_i(d\eta_{i}|\eta_{i^c})-\g'_i(d\eta_{i}|\bar\eta_{i^c})
\Big\Vert\leq \sum_{j\in G\ba i} Q_{i, j }d'(\eta_j,\bar\eta_j).
\end{split}
\end{equation}
where
\begin{equation}\begin{split}\label{dino}
& Q_{i, j }= 4  e^{2 \sum_{A\ni i}\Vert \P_A\Vert_\infty }\Bigl(
\sum_{k\in G\ba i}\d_k\Bigl( \sum_{A\supset \{i,k\}}\Phi_{A} \Bigr)
\bar D_{k j}\Bigr) e^{\sum_{A\ni j}\d_j(\P_A))}
\end{split}
\end{equation}
with $\bar D=\sum_{n=0}^\infty \bar C^n$.
\end{enumerate}
\end{thm}

Note that the first layer system may be very well in a
phase transition regime. For arbitrarily large interactions $\Phi$,
good concentration of the conditional measures $\a_{\eta_i}$ can
still lead to a small "Dobrushin matrix" $\bar C$, when the
deviation matrix $\bar B(\eta)$ is uniformly small in $\eta$. In
short: Uniform conditional Dobrushin uniqueness of the first layer
implies Gibbsianness of the second layer, with explicit estimates.

The matrix $Q$ describing the spatial loss of memory of the
variation of the conditioning, depends on the summability properties
of the potential $\Phi$ and the decay of the "Dobrushin-matrix"
$\bar C$. Note that the summability property we impose on the
initial potential (\ref{triple}) implies the finiteness of
(\ref{dino}).  In particular we have the following bound on the
entries of the $Q$-matrix;
\begin{equation}\label{qij}
Q_{ij}\leq 4\exp\Big(4\sup_{i\in G}\sum_{A\ni
i}||\P_A||_\infty\Big)\big(M\bar D\big)_{ij},
\end{equation}
where $M$ is the matrix  given by
$M_{ik}=\left\{\begin{array}{ll}
\d_k\big(\sum_{A\supset\{i,k\}}\P_A\big) &\mbox{if $ i\neq k$};\\
0  &\mbox{if $i= k$}.\end{array} \right. $

All these quantities are easily made explicit in examples.

\subsection{Goodness of short-time Gibbsianness for time-evolved rotator
 models}\label{Short-Time}

Let us get more concrete. 
Consider the rotator model on $G$, with both first
layer and second layer local spin spaces equal to $S^{q-1}$,  the
sphere in q-dimensional Euclidean space, with $q\geq2$. 

Take as a
Hamiltonian of the first layer system in infinite volume
\begin{equation}\label{rotor}
H(\s)=-\sum_{i,j\in G}J_{ij}\s_i\cdot\s_j
\end{equation}
with $\sup_{i}\sum_{j}|J_{ij}|<\infty$ where we assume that
$J_{ii}=0$ for
 each $i\in G$. Let $K$  be given by 
$K(d\s_i,d\eta_i)=K_t(d\s_i,d\eta_i)=k_t(\s_i,\eta_i)\a_0(d\s_i)\a_0(d\eta_i)$,
where $\a_0$ is the equidistribution on  $S^{q-1}$ and
$k_t$ is the heat kernel on the sphere, i.e.
\begin{eqnarray}
\Big(e^{\D t}\varphi\Big)(\eta_i)=\int \a_o(d\s_i)k_t(\s_i,\eta_i)\varphi(\s_i),
\end{eqnarray} 
where $\D$ is the Laplace-Beltrami operator on the sphere and $\varphi$ is any test function.  
$k_t$ is also called the \textbf{ Gauss-Weierstrass kernel}. The time-evolved
measure is given by
\begin{equation}\label{eq:intro.2}
\begin{split}
\mu_{t}(d\eta)=\int\mu(d\s)\prod_{i}k_t(\s_i,\eta_i)\a_0(d\eta_i).
\end{split}
\end{equation}
It has the product over the equidistributions on the spheres as an
infinite-time local limiting measure
\begin{equation}
\lim_{t\uparrow\infty}\mu_{t}(d\eta)=\bigotimes_{i\in G}
\a_0(d\eta_i).
\end{equation}
Denote $\g'_t$ by the class of all finite-volume conditional distributions 
of the time-evolved measure with full $\eta$-conditioning.
Then the following continuity estimates on the conditional probabilities
of the time-evolved model hold.

\begin{thm}  \label{thm:main-intro} Denote by $d(\eta,\eta')$ the
 induced metric
on the sphere $S^{q-1}$ (with $q\geq 2$) obtained by embedding the
sphere into the Euclidean space $\R^q$.

Assume that
\begin{equation}\begin{split}\label{rotatorc}
&\sqrt2\Bigl( \sup_{i } \sum_{j\in G } e^{ |J_{ij}|}|J_{ij}|\Bigr)
\left(1-e^{-(q-1)t} \right)^{\frac{1}{2}} <1.\cr
\end{split}
\end{equation}
Then the following holds.  
\begin{enumerate}
\item  The measure $\mu_t$ is Gibbs for a specification $\g'_t$, and

\item $\g_t$ satisfies the continuity estimate
\begin{equation}\begin{split}
&\Big\Vert\g'_{i,t}(d\eta_{i}|\eta_{i^c})-\g'_{i,t}(d\eta_{i}|\bar\eta_{i^c})\Big\Vert\leq
\sum_{j\in G\ba i} \bar Q_{i, j}(t)d(\eta_j,\bar\eta_j),
\end{split}
\end{equation}
with
\begin{equation}\begin{split}\label{barq}
& \bar Q_{i, j}(t):= \frac{1}{2}\min\Bigl\{ 
\sqrt{\frac{\pi}{t}}
Q_{i, j}(t), e^{4\sum_{l}|J_{jl}|}-1\Bigr \}
\end{split}
\end{equation}

where

\begin{equation}\begin{split}
& Q_{i, j }(t)= 8 e^{4\sup_{i\in G}
 \sum_{j\in G}|J_{ij}|}\sum_{k\in G\ba
 i}|J_{ik}|\bar D_{kj}(t),
\end{split}
\end{equation}
$\bar D(t) = \mathbf{1}+\sum_{n=1}^\infty \big(
1-e^{-(q-1)t}\big)^{\frac{n}{2}} A^n $,  $A$ is the matrix whose
entries are given by $A_{ij}=e^{|J_{ij}|}|J_{ij}|$  and $\mathbf{1}$
is the identity matrix.

\end{enumerate}
\end{thm}

In the definition of $Q_{ij}$ in the above theorem we have used
the bound (\ref{qij}) on the $Q_{ij}$ in Theorem \ref{main1}.

The proof of the theorem follows from three ingredients: 1) Theorem \ref{main1} 
which gives a continuity estimate in terms of the posterior metric $d'$, 
2) a comparison result between $d'$ and $d$, see Proposition \ref{Festimates}
and 3) a telescoping argument over sites in the conditioning. 

It is straightforward to apply Theorem \ref{main1} to our model and obtain 
a result formulated in $d'$. However, a more natural metric we would prefer 
to use is $d$, and so 
we should use a comparison argument, applying   Proposition \ref{Festimates}. 
What continuity estimates do we expect to gain from this? 
It is elementary to see that for the initial kernel
\begin{equation}\begin{split}
&\Big\Vert\g_{t=0}(d\eta_{i}|\eta_{i^c})-\g_{t=0}(d\eta_{i}|\bar\eta_{i^c})\Big\Vert\leq
 e^{2\sum_{j\in G}|J_{ij}|}\sum_{j\in G}|J_{ij}|d(\eta_j,\bar\eta_j)
\end{split}
\end{equation}
We see that continuity can be measured 
in terms of $d$, due to the Lipschitz property 
of the initial Hamiltonian,  and the spatial decay is provided by 
the decay of the couplings. 

So, at small  time $t$, we are aiming at a similar continuity estimate to
hold which is uniform in $t$ as $t$ goes to zero. 
Now, while estimating $d'$ against $d$ we have accumulated a nasty factor $\frac{1}{\sqrt{t}}$ that 
blows up when time $t$ goes to zero. We note that this is not just an artefact of 
Proposition \ref{Festimates}, but the posterior metric between two points
on the sphere indeed blows up like $\frac{1}{\sqrt{t}}$, as  can be seen from the proof. 
At first sight this does not seem to be a problem in the definition of $Q_{ij}(t)$ 
because the off-diagonal entries of the matrix $\bar D_{ij}(t)$ are 
suppressed by the same factor proportional 
to $\sqrt{t}$ that appears in (16).  This suppression follows from a bound 
on the corresponding Dobrushin matrix of this order. 
Unfortunately the diagonal terms of $\bar D(t)$ give rise to blow-up for sites $i$ and $j$ 
that are within the range of the potential. 
As it is clear from the proof, this blow-up is understandable since  
so far we did not employ any continuity properties of the initial
Hamiltonian w.r.t.  the Euclidean metric. Without further
conditions of this sort clearly no continuity can be expected, as
even a system of two sites with the Hamiltonian being a step
function shows. 

Now, to disentangle these {\it local effects} from the {\it global
effects} treated so far,  we use in the third step a telescoping argument over 
the conditioning. Exploiting Lipschitz-continuity w.r.t. a single argument of the Hamiltonian 
we obtain the second term in the minimum in  (\ref{barq}) which puts a time-independent ceiling
to the blow-up for small times. This solves the blow-up problem.

In this context let us also exhibit the comparison estimate of the two metrics $d'$ and $d$ that
we also deem of interest in itself.
\begin{prop} \label{Festimates}
There is an estimate of the posterior metric $d'$ associated to the
measure $K_t$ of the form
\begin{equation}\begin{split}
&d'(\eta_j,\bar\eta_j)\leq F_{q,t}\bigl (
 d(\eta_j,\bar\eta_j) \bigr).
 \end{split}
 \end{equation}

The function $F_{q,t}$ satisfies the following:

\begin{enumerate}
\item For any $q\geq 2$,  $x\in [0,2]$ and $t>0$ we have the estimate
\begin{equation}\begin{split}
F_{q,t}(x)&\leq 4 P\Bigl( 0\leq G\leq \frac{\arcsin
\frac{x}{2}}{\sqrt{2t}} \Bigr)\leq \frac{\sqrt{\pi}x}{2\sqrt{ t}}
\end{split}
\end{equation}
where $G$ is a standard normal variable.

\item In general dimensions $q\geq 2$ more information can be derived
 by
the expansion
\begin{equation}\begin{split}
F_{q,t}(x)&=\sum_{m=0}^\infty
a_{q,m}(t)P_{2m+1}\Big(q,\frac{x}{2}\Big)\text{
 with }
\cr a_{q,m}(t)&=e^{-(2m+1) (2m+q-1)t}\dfrac{(-1)^{m} 4N(q,m)\Gamma\left( \frac{q}{2}\right)
}{\sqrt{\pi}\Gamma\left( \frac{q-1}{2}\right)}
\prod_{i=0}^m \Bigg( \dfrac{2i-1}{ q+2i-1}\Bigg)
\end{split}
\end{equation}
in terms of  Legendre polynomials $P_n(q,s)$ of degree $n$ in
dimension $q$ (see the definition \ref{Rod})and $N(q,m)$ is also the dimension of the space of spherical harmonics of
degree $n$ in dimension $q$ (see (\ref{kt})).
\end{enumerate}
\end{prop}
{\bf Remark: } 
The proof uses a coupling argument and a reflection principle for
diffusions on the sphere under reflection at the equator.
\bigskip

\subsection{Goodness of Gibbsianness for local approximations}

As another consequence from the general theorem we prove 
that any sufficiently fine local coarse graining preserves the Gibbs 
property.  Here the fineness of the coarse graining 
has to be compared relative to the scale in the local state spaces on which 
the initial Hamiltonian is varying.  

We thus need a bit more structure, namely let $(S,d)$ now be a metric space. 
Let a decomposition be given of the form $S=\bigcup_{s'\in S'}S_{s'}$. Here 
$S'$ may be a finite or infinite set. 
Put $T(s):=s'$ for $S_{s'}\ni s$. This defines a deterministic 
transformation on $S$, called the {\it fuzzy map}.  With this map we associate to 
each $s'\in S'$ a corresponding a priori measure on $S_{s'}$ (say $\a_{s'}$ ). Note 
that here $\a_{s'}$ is the corresponding analogue of $\a_{\eta}=K(\cdot|\eta)$ for the fuzzy map.   

\begin{thm}\label{fuzzy}
Assume the Lipschitz-property for the $j$-variation of the initial Hamiltonian 
\begin{equation}\begin{split}\label{lippi}
& \sup_{{\zeta,\bar \zeta}\atop{\zeta_{j^c}=\bar \zeta_{j^c}}} \Bigl
| H_i(\s_i\z_{i^c})-H_i(\s_i\bar\z_{i^c})-\Big(H_i(a_i\z_{i^c})-
H_i(a_i\bar\z_{i^c})\Big)\Bigr |\leq L_{i j }d(\s_i,a_i).
\end{split}
\end{equation}
Suppose that 
\begin{eqnarray}
\frac{\rho}{2}\sup_{i\in G}\sum_{j\in G }\exp\Bigr(
\frac{1}{2}\sum_{A\supset\{i,j\}}\delta(\Phi_A)\Bigl)L_{ij}\nonumber < 1
\end{eqnarray}
where $\rho=\sup_{s'}\text{diam}(S_{s'})$ denotes the fineness of the decomposition. 

\begin{enumerate}
\item Then, for any initial Gibbs measure $\mu$ of the specification $\Phi$ with an 
arbitray a priori measure $\a$ the 
transformed measure $T(\mu)$ is Gibbs for a specification $\g'$. 

\item The entries  $C'_{ij}$ of the Dobrushin interdependence matrix of $\g'$
are bounded by $Q_{ij}$ given by (\ref{dino}) where we have to put 
\begin{eqnarray}
\bar C_{ij}= \frac{\rho}{2}\exp\Bigr(
\frac{1}{2}\sum_{A\supset\{i,j\}}\delta(\Phi_A)\Bigl)L_{ij}.\nonumber
\end{eqnarray}
\end{enumerate}

\end{thm}

Answering a question of Aernout van Enter, 
this provides a class of examples where $S$ and $S'$ are different (one may 
be continuous, the other not),  
the initial measure may be in the phase transition regime, and 
the image measure will be Gibbs. 
To think of an even more concrete example,
let take the rotor-model (\ref{rotor}). 
Divide the sphere $S^{q-1}=\bigcup_{s'}S_{s'}$ into "countries" $S_{s'}$.  
Then the correspondingly discretized 
model on the country-level is still Gibbs whenever there is no country 
with diameter bigger then 
$\Bigl( \sup_{i } \sum_{j\in G } e^{ |J_{ij}|}|J_{ij}|\Bigr)^{-1}$.

As a concluding remark let us mention that we may very well 
apply our method also to other well-known examples of transforms 
of Gibbs measures that may potentially lead 
to renormalization group pathologies. For instance, also 
the decimation transformation mapping a Gibbs measure on the lattice 
to its restriction to a sublattice can be cast in this framework. 
Theorem \ref{main1} then implies the statement that the projected measure is 
always Gibbs if the interaction is sufficiently small in triple norm. 
The posterior metric for configurations on the projected lattice 
then becomes the discrete metric 
$d'(\eta_i,\eta_i')=1_{\eta_i\neq \eta'_i}$ and hence 
the matrix element $Q_{ij}$ becomes a bound on the Dobrushin interdependence matrix 
of the image system.

\bigskip 
\bigskip
\section{On the proofs on Theorem \ref{Dob3} and \ref{concen}:}
In this section we provide proofs of Theorem \ref{Dob3} and
\ref{concen} and also state and prove some related results.
We start with the \\\\
\textbf{Proof of Theorem \ref{Dob3}:} The idea of the proof is to
find an estimate on  the Dobrushin interdependence matrix  as in the
proof of Proposition 8.8  of \cite{GOR}. This involves  estimating
the variation of the single-site measure  at a given site $i\in G$
when varying the boundary condition at some site $j\in G\backslash
i$. That is we fix $\z,\eta\in \O$ with $\z_{j^c}=\eta_{j^c}$  and
put $u_0(\s_i)=-H_i(\s_i\zeta_{i^c})$ and
$u_1(\s_i)=-H_i(\s_i\eta_{i^c})$. We proceed further by taking the
linear interpolation $u_t=t u_1+(1-t)u_0$ of $u_1$ and $u_0$. It
follow from this linear interpolation  that
 \begin{equation}\begin{split}\label{ut}
\d(u_t)\leq \sum_{A\supset\{i,j\}}\d(\Phi_A).
\end{split}
\end{equation}
Setting  $h_t=e^{u_t}/\a(e^{u_t})$ and
$\l_t(d\s_i)=h_t(\s_i)\a(d\s_i)$ we note that
$\l_0(d\s_i)=\g_i(d\s_i|\z)$ and
 $\l_1(d\s_i)=\g_i(d\s_i|\eta)$.
We now observe that
 \begin{eqnarray}\label{linear}
 2\Vert\l_0-\l_1 \Vert =\int \a(d\s_i) |h_1(\s_i)-h_0(\s_i)|= \int
 \a(d\s_i) \Bigl |\int_{0}^1 dt\frac{d}{dt}h_t(\s_i)\Bigr |\nonumber\\
\leq \int_0^1 dt\;
\l_t\Bigr|H_i(\cdot\z_{i^c})-H_i(\cdot\eta_{i^c})-\l_t\Big(H_i(\cdot\z_{i^c})-H_i(\cdot\eta_{i^c})\Big)\Bigl|
\\
\leq 2 \int_0^1 dt \exp\Bigr(\d(u_t)\Bigl)  \inf_{B}\int \a(d\s_i)
\Bigr|H_i(\s_i\z_{i^c})-H_i(\s_i\eta_{i^c})-B\Bigl|.\nonumber
\end{eqnarray}
It  follows from(\ref{Ci}),(\ref{deval}) and (\ref{ut})  that
\begin{equation*}\begin{split}
& C_{ij} \leq \exp\Bigr(\sum_{A\supset\{i,j\}}\d(\Phi_A)\Bigl)
\text{dev}_{\a;i,j}(H_i).
\end{split}
\end{equation*}
The rest of the proof follows from the definition of the Dobrushin
constant $c$.
\begin{flushright}
$\Cox$
\end{flushright}
Sometimes  it is useful to  use  quadratic variation instead of  the
linear variation $ \text{dev}_{\a;i,j}$ to obtain an explicit bound,
as we shall  see in the proof of Theorem \ref{thm:main-intro} below.
We define this quadratic variation as follows.
\begin{defn}\label{quad}
For any bounded measurable function $F$ on $\O$ we define for any
pair $i,j\in G$  $\text{std}_{\a;i,j}(F)$ as
\begin{eqnarray}
\text{std}_{\a;i,j}(F):=\sup_{\z,\bar\z\in \O,
\z_{j^c}=\bar\z_{j^c}} \inf_B\Bigg(\int d\a(d\s_i)\Big(
F(\s_i\z_i^c)-F(\s_i\bar\z_i^c)-B\Big)^2\Bigg)^{\frac{1}{2}}.
\end{eqnarray}
\end{defn}
The quantity  $\text{std}_{\a;i,j}(F)$ is the worst-case quadratic
deviation of the variation of $F$ at the site $j$ viewed as a random
variable w.r.t. to $\s_i$ under $\a(d\s_i)$. Clearly
$\text{dev}_{\a;i,j}(F)\leq\text{std}_{\a;i,j}(F)$, so we could
bound the inequality in Theorem \ref{Dob3} in terms of the quadratic
variation; going directly into the proof however gives a slightly
better constant.

This gives rise to the following "quadratic version" of Theorem
\ref{Dob3}.
\begin{prop}\label{corBob3}
The Dobrushin constant $c$ is also bounded by
\begin{equation}\begin{split}
c\leq \frac{1}{2}\sup_{i\in G}\sum_{j\in G}\exp\Bigr(\frac{1}{2}
\sum_{A\supset\{i,j\}}\delta(\Phi_A)\Bigl)\text{std}_{\a;i,j}(H_{i}).\cr
\end{split}
\end{equation}
\end{prop}
\textbf{Proof:} The proof uses the same arguments employed in the
proof of Theorem \ref{Dob3} above, the only difference being that we
have a quadratic estimate  (resulting from the Cauchy-Schwartz inequality)  in
\begin{eqnarray}\label{cstd}
 2\Vert\l_0-\l_1 \Vert =\int \a(d\s_i) |h_1(\s_i)-h_0(\s_i)|= \int
 \a(d\s_i) \Bigl |\int_{0}^1 dt\frac{d}{dt}h_t(\s_i)\Bigr |\nonumber\\
\leq \int_0^1 dt\;
\l_t\Bigr|H_i(\cdot\z_{i^c})-H_i(\cdot\eta_{i^c})-\l_t\Big(H_i(\cdot\z_{i^c})-H_i(\cdot\eta_{i^c})\Big)\Bigl|
\\
\leq \int_0^1 dt \exp\Bigr(\frac{\d(u_t)}{2}\Bigl)\inf_{B}\Bigr(
\int\a(d\s_i)\Bigr(H_i(\s_i\z_{i^c})-H_i(\s_i\eta_{i^c})-B\Bigl)^2\Bigl)^{\frac{1}{2}}.
\nonumber
\end{eqnarray}
\begin{flushright}
$\Cox$
\end{flushright} 

If the initial Hamiltonian satisfies a Lipschitz-property w.r.t. 
 a given metric $d$ on the local state space an estimate 
 of the Dobrushin constant can be formulated as follows.

\begin{cor} \label{Lips}
Suppose the Lipschitz-condition (\ref{lippi}).
Then we have
\begin{equation*}
 c\leq \frac{1}{2}\sup_{i\in G}\sum_{j\in G}\exp\Bigr(\frac{1}{2}
\sum_{A\supset\{i,j\}}\delta(\Phi_A)\Bigl)L_{ij}\inf_{a_i\in
S}\Big(\int d^2(\s_i,a_i)\a(d\s_i) \Big)^{\frac{1}{2}}.
\end{equation*}
\end{cor}
The proof of the corollary follows from the last inequality in
(\ref{cstd}), since taking the infimum over $B$ is less than or equal to
taking the infimum over $a_i\in S$ when we substitute $B$ with $H_i(a_i\z_{i^c})-
H_i(a_i\eta_{i^c})$.

A somewhat more abstract reformulation of the bounds on the
Dobrushin's constant  can be given in terms of appropriately defined
norms of the potential.
\begin{cor}
Define for  $\P$ the norms
\begin{equation}\begin{split}\label{tippnorm}
&\triplenorm\P\triplenorm_{dev_{\a}}:=\sup_{i\in G}\sum_{j\in
G\setminus \{i\}}\sum_{A\supset\{i,j\}}\text{dev}_{\a;i,j}\left(
\P_A\right)\cr &\triplenorm\P\triplenorm_{std_{\a}}:=\sup_{i\in
G}\sum_{j\in G\setminus
\{i\}}\sum_{A\supset\{i,j\}}\text{std}_{\a;i,j}\left( \P_A\right).
\end{split}
\end{equation}
Then we have for $\triplenorm\P\triplenorm<\infty$ that the
Dobrushin constant $c$ of the specification for $\P$ satisfies the following bounds
 \begin{eqnarray}
c\leq
e^{2\triplenorm\P\triplenorm}\triplenorm\P\triplenorm_{dev_{\a}}\quad
\text{and}\quad c\leq
\frac{1}{2}e^{\triplenorm\P\triplenorm}\triplenorm\P\triplenorm_{Std_{\a}}.
\end{eqnarray}
\end{cor}
Once the definitions are made the proof is obvious. Note further
that $c$ is  finite as long as $\triplenorm\P\triplenorm$ is because
$\triplenorm\P\triplenorm_{std_{\a}}\leq 2
\triplenorm\P\triplenorm$. We finally give the proof of the bounds
in the "Concentration implies concentration"-theorem.\\\\
\textbf{Proof of Theorem \ref{concen}:} Note 
that the hypothesis $||B||_\infty<\frac{1}{s}$
(as we will see below)  implies that we are in the uniqueness regime.
Then for any bounded measurable
function $F$ on $\O$  it follows from Theorem 1 of \cite{KUL3} that
under the unique Gibbs measure $\mu$
\begin{eqnarray}
\mu\Big(F-\mu(F)\geq r\Big)\leq
\exp\Bigg(-\dfrac{r^2}{2}\dfrac{(1-c)(1-c_t)}{||\underline{\d}(F)||^2_{l_2}}
\Bigg) \quad \forall\; r\geq0,
\end{eqnarray}
where $c$ and $c_t$ are respectively the Dobrushin constants of
the Dobrushin interdependence matrix and its transpose.
It follows from the definitions of $c$ and $c_t$ and the bound in 
Theorem \ref{Dob3} that
\begin{eqnarray}
c &\leq& \exp\Bigr( \sup_{i\neq
 j}\sum_{A\supset\{i,j\}}\delta(\Phi_A)\Bigl)\sup_{i\in G}\sum_{j\in
 G}\text{dev}_{\a;i,j}(H_{i})=s ||B||_\infty\nonumber\\
c_t &\leq & \exp\Bigr( \sup_{i\neq
j}\sum_{A\supset\{i,j\}}\delta(\Phi_A)\Bigl)\sup_{j\in G}\sum_{i\in
G}\text{dev}_{\a;i,j}(H_{i})=s ||B||_1\nonumber.
\end{eqnarray}
\begin{flushright}
$\Cox$
\end{flushright}
Note that the validity of Theorem \ref{concen} depends on $c$ and
$c_t$ being less than one. In our criterion (\ref{c-bound}) the
smallness of the $\text{dev}_{\a;i,j}$'s is the main  ingredient for
$c$ and $c_t$ to be less than one. This smallness of the
$\text{dev}_{\a;i,j}$'s is caused by good "concentration" properties
of $\a$ even if the  interaction is strong and possibly by the weakness 
of the interaction.

\section{ On the proof Theorem \ref{main1} and related results} The 
purpose of
this section is to give the proof  of Theorem \ref{main1}  outlined
in Section \ref{twolayer}  of the introduction. The main ingredient
to the proof is to show the  lack  of phase transitions in some intermediate
system and exploit the consequences for decay of spatial memory. 
Recall from Section \ref{twolayer} that our  initial system
was given by the Gibbs measure $\mu$ admitted by the specification
$\g$ obtained from the
 interaction $\Phi$ and an a priori measure $\a=\int K(\cdot,d\eta_i)$
described above.  Thus for a given boundary condition $\bar \s$ and
any finite volume $\L\subset G$  we write $\g_\L(\cdot|\bar\s)\in
\g$ as
\begin{eqnarray}
\g_\L(d\s_\L|\bar\s)= \frac{ \exp\Bigl( - H_{\L}(\s_\L\bar
\s_{\L^c}) \Bigr) \prod_{j\in \L}\a(d\sigma_j) }{\int_{S^\L}
\exp\Bigl( - H_{\L}(\td\s_\L\bar \s_{\L^c}) \Bigr) \prod_{j\in
\L}\a(d\td\s_j )} ,
\end{eqnarray}

We now introduce a double-layer system or  joint system  by coupling
the initial system  to a second system (with single-spin space $S'$)
through the sitewise joint  measures $K(d\s_i,d\eta_i)$ on $ S\times
S'$. Denote by $\td\g$ the specification of our new double-layer
system, i.e. for a fixed boundary condition $\bar \s\in\O=S^G$ and a
finite volume $\L\subset G$, $\td\g_\L(\cdot|\bar \s)\in\td\g$ is
given by

\begin{eqnarray}\label{gamatd}
\td\g_\L(d\s_\L,d\eta_\L|\bar \s)&=& \frac{  \exp\Bigl( -
H_{\L}(\s_\L\bar \s_{\L^c}) \Bigr) \prod_{j\in
\L}K(d\sigma_j,d\eta_j) }{\int_{(S\times S')^\L} \exp\Bigl( -
H_{\L}(\td\s_\L\bar \s_{\L^c})
\Bigr) \prod_{j\in \L}K(d\td\s_j,d                \td\eta_j )}\\
 &=& \g_\L(d\s_\L|\bar\s)\prod_{j\in \L}K(d\eta_j|\s_j),\nonumber
\end{eqnarray}
where $K(d\eta_i |\s_i)$ denotes  the $K$ conditional distribution
of the second spin given the value of the first. This specification
is in general not Gibbs but in our case where we only have sitewise
dependence between the two layers it is known for instance from
\cite{ACD} and references therein that $\td\g$ is Gibbs.

For each non-empty  subset $\L$ of $G$ we denote by $\mathcal{S}_\L$
the collection of all non-empty finite subsets of $\L$. We will
write $\mathcal{S}$ instead of $\mathcal{S}_G$. For any fixed
configuration $\s\in \O=S^G$ and any $\L\in \mathcal{S}$ we define
the finite-volume \textit{transformed } distribution
$\g'_{\L,\bar\s}$ as
\begin{eqnarray}\
\g'_{\L;\bar\s}(d \eta_\L)&:=& \int_{S^\L}\td\g(d\s_\L,d\eta_\L|\bar
\s_{\L^c}) .
\end{eqnarray}
It is important to note that in the joint system considered above,
conditionally on the $\s$'s the $\eta$'s are independent. But taking the
$\s$-average of the joint system creates dependence among the
$\eta$'s. Due to this dependence we now introduce finite-volume
$\eta$ conditional distributions by freezing  the  $\eta$
configuration in the definition of $\g'_{\L;\bar\s}$ except at some
region $\Delta\in\mathcal{S}_\L$. That is for any $\L\in\mathcal{S}$
with $|\L|\geq 2$ and $\Delta\in\mathcal{S}_\L$  we have
\begin{eqnarray}\label{1-pt}
\g'_{\D,\L; \bar \s}(d\eta_\D|\bar \eta_{\L\ba \D})= \frac{ \int_{S^\L}
\exp\Bigl( - H_{\L}(\s_\L\bar \s_{\L^c}) \Bigr) \prod_{j\in \L\ba
\D}K(d\sigma_j |\bar\eta_j)\prod_{i\in\D} K(d\sigma_i, d\eta_i)
}{\int_{S^\L} \exp\Bigl( - H_{\L}(\s_\L\bar \s_{\L^c}) \Bigr)
\prod_{j\in \L\ba \D}K(d\sigma_j |\bar\eta_j)\prod_{i\in\D}
\a(d\sigma_i ) }.\cr
\end{eqnarray}
The natural question that comes to mind is whether $\lim_{\L\uparrow
G}\g'_{\D,\L; \bar \s}(d\eta_\D| \bar\eta_{\L\ba \D})$ exists for any
fixed $\D\in \mathcal{S}_\L$, $\bar\s\in \O$  and
$\bar\eta_{\D^c}\in (S')^{G\ba \D}$? If this limit exists we will
denote it by $\g'_\D(d\eta_\D|\bar\eta_{\D^c})$ and $\g'$ by the class 
of all the conditional distributions for finite $\D$. For the sake of
simplicity we will always restrict our analysis to the case where
$\D$ is a singleton. The analysis for general (but finite) $\D$ can
be implemented using the same arguments used in the singleton case.
It is our aim to provide a sufficient condition for the conditional
probabilities $\g'_{i,\L; \bar \s}(d\eta_i| \eta_{\L\ba i})$ to have
an infinite-volume limit. For this we introduce the decomposition of
the Hamiltonian $H_{\L}$ in the finite window $\L$ into its
contributions coming from the sites in $\L\ba i$ and  site $i$ for
any $i\in \L$ as follows;
\begin{eqnarray}\label{hdecomp}
H_{\L}(\s_{\L}\bar\s_{\L^c})&=& H_{i}(\s_{\L}\bar\s_{\L^c})+
H_{\L\ba i}(\s_{\L\ba i}\bar\s_{\L^c}),\quad \text{where}\nonumber\\
 H_{i}(\s_{\L}\bar\s_{\L^c})&=&\sum_{A\ni i}\P_A(\s_{\L}\bar\s_{\L^c})
 \quad\text{and}\\
 H_{\L\ba i}(\s_{\L\ba i}\bar\s_{\L^c})&=&\sum_{A\cap \L\setminus i
 \neq\emptyset ;\; i\notin A }\P_A(\s_{\L\ba i}\bar\s_{\L^c})\nonumber
\end{eqnarray}
 We clearly see from the
definition of an interaction that the Hamiltonian $ H_{\L\ba i}$ is
a function on the configuration space $S^{G\setminus i}$. For the
infinite-volume transformed conditional distributions $\g'_i(d\eta_i|
\eta_{ i^c})$ to exist, some intermediate system living on the
sublattice $G\ba i$  must admit a unique infinite-volume Gibbs
measure. This intermediate model is what we referred to as the
\textit{\textbf{restricted constrained first layer model}} (defined
below w.r.t $ H_{\L\setminus i}$ ).

\begin{defn}
The \textit{\textbf{restricted constrained first layer model}}
(\textbf{RCFLM}) in any  $\L\in \mathcal{S}$ with $|\L|\geq 2$ and
$i\in\L$ is defined as the measure,
\begin{eqnarray}
\mu^{\bar\s}_{\L\ba i}[\eta_{\L\ba i}](d\s_{\L\ba i})= \frac{
\exp\Bigl( - H_{\L\ba i}(\s_{\L\ba i}\bar\s_{\L^c}) \Bigr)
\prod_{j\in \L\ba i}K(d\sigma_j |\eta_j) }{\int_{S^{\L\ba i}}
\exp\Bigl( - H_{\L\ba i}(\tilde \s_{\L\ba i}\bar \s_{\L^c}) \Bigr)
\prod_{j\in \L\ba i}K(d\tilde \sigma_j |\eta_j) },
\end{eqnarray}
for some $\bar\s=S^{G}$ and $\eta_{\L}\in (S')^{\L}$.
\end{defn}
It is \textbf{restricted} because we only consider the spins in the
sublattice $G\setminus i$ and \textbf{constrained} since we have
frozen the configuration in the second layer . The RCFLM (as we will
see from the lemma below) will provide us with a sufficient
condition for the existence of an infinite-volume limit
$\g'_i(d\eta_i| \eta_{i^c})$ for the conditional probabilities
$\g'_{i,\L; \bar \s}(d\eta_i| \eta_{\L\ba i})$.

\begin{lem}\label{kenel}
Let $\L\in \mathcal{S}$ with $|\L|\geq2$, then for any $i\in \L$ and
any $\bar\s\in \O$ we have
\begin{eqnarray}
\g'_{i,\L; \bar \s}(d\eta_i| \eta_{\L\ba i})= \dfrac{\int_{S^{\L\ba
i}} \mu^{\bar\s}_{\L\ba i}[\eta_{\L\ba i}](d\s_{\L\ba i})
\int_{S}\exp\Bigr(-H_{i}(\s_\L\bar\s_{\L^c})\Bigr) K(
d\s_i,d\eta_i)}{\int_{S^{\L\ba i}} \mu^{\bar\s}_{\L\ba
i}[\eta_{\L\ba i}](d\s_{\L\ba i})
\int_{S}\exp\Bigr(-H_{i}(\s_{\L}\bar\s_{\L^c})\Bigr) d\a( \s_i)} .
\end{eqnarray}
\end{lem}
\textbf{Proof:} By using the decomposition of $H_{\L}$ in
(\ref{hdecomp}) we can write  $\g'_{i,\L; \bar \s}(d\eta_i|
\eta_{\L\ba  i})$ as;
\begin{equation}\begin{split}
&\g'_{i,\L; \bar \s}(d\eta_i| \eta_{\L\ba  i})= \cr
&\dfrac{\int_{S^{\L\ba i}} \exp\Bigl( - H_{\L\ba i}(\s_{\L\ba i}\bar
\s_{\L^c}) \Bigr) \prod_{j\in \L\ba i}K(d\sigma_j |\eta_j)
\int_{S}\exp\Bigr(-H_{i}(\s_{\L}\bar\s_{\L^c})\Bigr) K(
d\s_i,d\eta_i)}{\int_{S^{\L\ba i}} \exp\Bigl( - H_{\L\ba
i}(\s_{\L\ba i}\bar \s_{\L^c}) \Bigr) \prod_{j\in \L\ba
i}K(d\sigma_j |\eta_j) \int_{S\times
S'}\exp\Bigr(-H_{i}(\s_{\L}\bar\s_{\L^c})\Bigr) K(d \s_i,
d\tilde{\eta_i})}.
\end{split}
\end{equation}
The claim of the lemma follows by multiplying the expression for
$\g'_{i,\L; \bar \s}(d\eta_i| \eta_{\L\ba i})$  above by
$\dfrac{\int_{S^{\L\ba i}} \exp\Bigl( - H_{\L\ba i}(\tilde \s_{\L\ba
i}\bar \s_{\L^c}) \Bigr) \prod_{j\in \L\ba i}K(d\tilde \sigma_j
|\eta_j) }{\int_{S^{\L\ba i}} \exp\Bigl( - H_{\L\ba i}(\tilde
\s_\L\bar \s_{\L^c}) \Bigr) \prod_{j\in \L\ba i}K(d\tilde \sigma_j
|\eta_j) }$ and simplifying the resulting expression.
\begin{flushright}
$\Cox$
\end{flushright}

It is not hard to infer from the above lemma  that there will be an
infinite-volume kernel $\g'_i(d\eta_i|\eta_{i^c})$ if the RCFLM has a
unique infinite-volume Gibbs measure $\mu_{i^c}[\eta_{i^c}]$. This
is the case since $H_i$ is a local function which is finite by
assumption. This was also observed in the corresponding mean-field
set-up  in \cite{CA}. Over there a sufficient condition for the
existence of infinite-volume transformed kernel was given in terms
of the uniqueness of global minimizers for some potential function.
This condition was shown to be equivalent to the differentiability
of the transformed Hamiltonian. We now  state a result concerning an
upper bound for  Dobrushin's constant for the RCFLM.

\begin{prop}\label{2.2}
Let the Dobrushin's interdependence  matrix for the RCFLM for some
fixed site $i_o\in G$ be the  matrix whose entries are given by
 \begin{eqnarray}
C^{i_o}_{ij}[\eta_i]=\sup_{\z,\bar\z\in S^{G\sm
i_o};\;\zeta_{j^c}=\bar\z_{j^c}}\Big\Vert  \mu^{ \z}_{i}[\eta_{i}] -
\mu^{\bar\z}_{i}[\eta_{i}] \Big\Vert,
\end{eqnarray} for any pair $i,j\in G\ba i_o$ where we have
denoted $ \mu^{\bar\z}_{i}$ by the single-site  part of $\mu^{\bar
 \z}_{\L\ba{i_o}}$ .\\
Then we have;
\begin{equation}\begin{split}
&C^{i_o}_{ij}[\eta_i] \leq\exp\Bigr(
\sum_{A\supset\{i,j\};\;i_o\notin
A}\delta(\Phi_A)\Bigl)\text{dev}_{\a_{\eta_i};i,j}(H_i).\cr
\end{split}
\end{equation}
where   $\a_{\eta_i}(d\s_i)=K(d\s_i|\eta_i)$. \\
Furthermore, defining  the  Dobrushin constant $c'[\eta]$ for the
RCFLM as
\begin{eqnarray}
c'[\eta]&:=&\sup_{i_o\in G}c^{i_o}[\eta], \quad \text{with}
\nonumber\\
c^{i_o}[\eta]&=&\sup_{i\in G\sm i_o}\sum_{j\in G \ba
i_o}C^{i_o}_{ij}[\eta_i],
\end{eqnarray}
 we also have
\begin{eqnarray}
c'[\eta]&\leq& \sup_{i_o\in G} \sup_{i\in G\sm i_o}\sum_{j\in G \ba
i_o}\exp\Bigr( \sum_{A\supset\{i,j\};\;i_o\notin
A}\delta(\Phi_A)\Bigl)\text{dev}_{\a_{\eta_i};i,j}(H_i)\cr &\leq &
\sup_{i\in G }\sum_{j\in G }\exp\Bigr(
\sum_{A\supset\{i,j\}}\delta(\Phi_A)\Bigl)\text{dev}_{\a_{\eta_i};i,j}(H_i)).
\end{eqnarray}
In the case of $G=\Z^d$ and translation-invariant initial interactions the
last inequality is an equality.

\end{prop}

\bigskip

{\bf Proof: } The proof follows the same lines as the proof of
Theorem \ref{Dob3} but here  we  use   $\a_{\eta_i}=K(\cdot|\eta_i)$
instead of $\a$ .
\begin{flushright}
$\Cox$
\end{flushright}

It is also not hard  to deduce  from Proposition \ref{corBob3} that;
\begin{eqnarray}
c'[\eta]\leq \frac{1}{2}\sup_{i_o\in G} \sup_{i\in G\sm
i_o}\sum_{j\in G \ba i_o}\exp\Bigr(
\frac{1}{2}\sum_{A\supset\{i,j\};\;i_o\notin
A}\delta(\Phi_A)\Bigl)\text{std}_{\a_{\eta_i};i,j}(H_i)\\\leq \frac{1}{2}\sup_{i\in G}\sum_{j\in G }\exp\Bigr(
\frac{1}{2}\sum_{A\supset\{i,j\}}\delta(\Phi_A)\Bigl)\text{std}_{\a_{\eta_i};i,j}(H_i)\nonumber.
\end{eqnarray}

Again Lipschitzness of the initial Hamiltonian carries over nicely.

\begin{cor}\label{maincor}Suppose the Lipschitz-condition (\ref{lippi}).
Then we have
\begin{equation}\begin{split}
 &c'[\eta]=\sup_{i_o\in G}\sup_{i\in G\ba i_o}\sum_{j\in
G\sm i_o}C^{i_o}_{ij}[\eta_i]\cr &\leq \frac{1}{2}\sup_{i_o\in
G}\sup_{i\in G\ba i_o}\sum_{j\in G\sm i_o}\exp\Bigr(
\frac{1}{2}\sum_{A\supset\{i,j\};\;i_o\notin
A}\delta(\Phi_A)\Bigl)L_{ij}\inf_{a_i\in S}\Bigl(\int_S
d^2(\s_i,a_i)\a_{\eta_i}(d\s_i)\Bigr)^{\frac{1}{2}}\cr
&\leq \frac{1}{2}\sup_{i\in G}\sum_{j\in G}\exp\Bigr(
\frac{1}{2}\sum_{A\supset\{i,j\}}\delta(\Phi_A)\Bigl)L_{ij}\inf_{a_i\in S}\Bigl(\int_S
d^2(\s_i,a_i)\a_{\eta_i}(d\s_i)\Bigr)^{\frac{1}{2}}.
\end{split}
\end{equation}
\end{cor}
The claim of the corollary  follows from Corollary \ref{Lips}.
\bigskip

We now proceed to prove  Theorem \ref{main1}, but before we do this
we still need  some results  from which the proof will follow. As a
first step we recall some known results about Dobrushin's uniqueness
concerning  an estimate of the distance between the unique Gibbs
measure admitted by a Gibbs specification satisfying  Dobrushin's
condition and  another Gibbs measure corresponding to some other
specification. This estimate tells us the local variation between the
two  infinite-volume probability measures.  This result which we state in the
proposition below can be found for example in \cite{GOR} as Theorem
8.20. Before we state the result we fix some notations. Suppose
$C(\g)$ is the Dobrushin interdependence matrix of a specification
$\g$ and $C^n(\g),\;n\geq 0$, the $n`$th power of $C(\g)$, then we
define the matrix
\begin{eqnarray}
D(\g)=(D_{ij})_{i,j\in G}:=\sum_{n\geq 0}C^n(\g). 
\end{eqnarray}

\begin{prop}\label{Dobdist}
Let $\g$ and $\bar\g$ be any two specifications with $\g$ satisfying
 Dobrushin's condition. Suppose that for each $i\in G$ we  have a
measurable function $b_i$ on the standard Borel space $\O$  with the
property that
\begin{eqnarray}\label{38}
||\g_i(\cdot|\s_{i^c})-\bar\g_i(\cdot|\s_{i^c})||\leq b_i(\s)
\end{eqnarray}
for all $\s\in \O$. Then for $\mu\in \GG(\g)$ and $\bar\mu\in
\GG(\bar\g)$ we have
\begin{eqnarray}\label{39}
|\mu(f)-\bar\mu(f)|\leq \sum_{i,j\in
G}\d_{i}(f)D_{ij}(\g)\bar\mu(b_j)
\end{eqnarray}
for all functions $f$ which are the uniform limits of functions that
depend on finitely many local  variables $\s_i$.
\end{prop}

Observe from Lemma \ref{kenel} that if the $RCFLM$  satisfies
Dobrushin's condition uniformly in $\eta$ the infinite-volume single-site 
kernels $\g'_i(\cdot|{\eta}_i^c)$ exist for every
$\eta$ . We will adapt the result in Proposition \ref{Dobdist} to
our present set-up to compare $\g'_i(\cdot|\eta_{i^c})$ and
$\g'_i(\cdot|\bar{\eta}_{i^c})$ for any pair of configurations $\eta,
\bar\eta\in \O'=(S')^G$. Further we denote  by $\g[\eta_{i^c}]$ the 
 specification of the RCFLM with full $\eta_{i^c}$ configuration. Again 
 we assume for the first layer model that  $\mu=\lim_n \mu_{\L_n}^{\bar\s}$
 as in the hypothesis of Theorem \ref{main1}.

\begin{prop}\label{37}
Suppose the RCFLM on the sublattice $G\ba i$ (for some $i\in G$)
satisfies Dobrushin's condition uniformly in $\eta$  with
unique infinite-volume limit $\mu_{i^c}[\eta_{i^c}]$. Then
\begin{enumerate}
\item the second layer system (the transformed
model) has infinite-volume single-site conditional distributions
$\g'_i(d\eta_i| \eta_{i^c})$ given by
\begin{eqnarray}\label{gamap}
\g'_i(d\eta_i| \eta_{i^c})=
 \dfrac{\int_{S^{G\ba i}} \mu_{i^c}[\eta_{i^c}](d\s_{i^c})
\int_{S}\exp\Bigr(-H_{i}(\s_i\s_{i^c})\Bigr) K(
d\s_i,d\eta_i)}{\int_{S^{G\ba i}} \mu_{i^c}[\eta_{i^c}](d\s_{i^c})
\int_{S}\exp\Bigr(-H_{i}(\s_i\s_{i^c})\Bigr) \a( d\s_i)}
\end{eqnarray}
\item
for any pair $\eta_{i^c},\bar\eta_{i^c}\in (S')^{G\ba i}$ we have
for any $j\neq i$ that
\begin{equation}\begin{split}
&\Big\Vert\g_j[\eta_{j}](\cdot|\bar\s_{G\ba i})-\g_j[\bar\eta_{j}](\cdot|\bar\s_{G\ba
i})\Big\Vert \leq 2\exp\Bigl(\sum_{A\ni
j}\d_j(\P_A)\Bigr)\Big\Vert
K(\cdot|\eta_j)-K(\cdot|\bar\eta_j)\Big\Vert,\cr
\end{split}
\end{equation}
where the $\g_j[\eta_{j}](\cdot|\bar\s_{G\ba i,})$'s are the single-site
 parts of the specification for the RCFLM for $i\in G$ and
$\eta_{i^c}$, and
\item given $h_2(\s_{i^c})=\int_{S\times
S'}K(d\s_{i},d\eta_{i})\exp\Bigl(H_{i}(\s_{i}\s_{i^c})\Bigr)$  it
follows that
\begin{equation}\begin{split}\label{deno3}
&\Big|\mu_{i^c}[\eta_{i^c}](h_2)-\mu_{i^c}[\bar\eta_{i^c}](h_2)\Big|
\leq 2 e^{\sum_{A\ni i}\Vert \P_A\Vert_\infty}\cr
&\times\sum_{k,j\in G\ba i}\d_k\Big(\sum_{A\supset \{i,k\}}\P_A\Big)
\Bar D_{kj} \exp\Bigl(\sum_{A\ni j}\d_j(\P_A)\Bigr)\Big\Vert
K(\cdot|\eta_j)-K(\cdot|\bar\eta_j)\Big\Vert.
\end{split}
\end{equation}
\item Furthermore, for any $k\neq i$   it is the case that

\begin{equation}\begin{split}
\d_k\Bigl( h_2(\s_{i^c})\Bigr) \leq \d_k\Bigl( \sum_{A\ni
i,k}\P_A\Bigr) e^{   \sum_{A\ni i}\Vert \P_A\Vert_\infty}
\end{split}
\end{equation}
\item  and  finally
\begin{equation}\begin{split} \label{last}
&\Big\Vert\g'_i(d\eta_{i}|\eta_{i^c})-\g'_i(d\eta_{i}|\bar\eta_{i^c})\Big\Vert
\leq
2\frac{\Big|\mu_{i^c}[\eta_{i^c}](h_2)-\mu_{i^c}[\bar\eta_{i^c}](h_2)\Big|}
{\mu_{i^c}[\bar\eta_{i^c}](h_2)}
\end{split}
\end{equation}
\end{enumerate}
\end{prop}
\textbf{Remark:} In particular, we can write for any finite volume  the corresponding relation 
 for the finite-volume conditional distribution with full $\eta$-conditioning as in (\ref{gamap}), i.e. if $\D\in \mathcal{S}$ the we have
\begin{eqnarray}\label{gamapg}
\g'_\D(d\eta_\D| \eta_{\D^c})=
 \dfrac{\int_{S^{G\ba \D}} \mu_{\D^c}[\eta_{\D^c}](d\s_{\D^c})
\int_{S^\D}\exp\Bigr(-H_{\D}(\s_\D\s_{\D^c})\Bigr)\prod_{i\in\D} K(
d\s_i,d\eta_i)}{\int_{S^{G\ba \D}} \mu_{\D^c}[\eta_{\D^c}](d\s_{\D^c})
\int_{S^\D}\exp\Bigr(-H_{\D}(\s_\D\s_{\D^c})\Bigr)\prod_{i\in \D} \a( d\s_i)}.
\end{eqnarray}

\textbf{Proof:}  \begin{enumerate}
\item The proof follows from a two-step limiting procedure. We fix an $\eta$-conditioning 
only in a finite volume $\D$ and construct the infinite-volume measure  of the RCFLM by fixing a boundary 
condition on the first layer outside $\L$ (which we assume for simplicity to contain $\D$) and let $\L$  tend to infinity. 
Then we let $\D$ tend to infinity, and recover the conditional probabilities by Martingale 
convergence and uniform approximation of the infinite-volume RCFLM, with conditionings only in volume $\D$. 

More precisely, it follows as in  Lemma \ref{kenel} that  we have for finite-volume 
conditionings the representation 

\begin{equation}
\g'_{i,\D,\L,\bar\s}(d\eta_{i}|\eta_{\D\ba i})
= \frac{\mu_{\L\ba i}^{\bar\s}[\eta_{\D\ba i}]\Big[\int_S e^{-H_i(\s_i\bar{\s}_{\L^c}\cdot_{\L\ba i})}K(d\s_i,d\eta_i)
\Big]}{\mu_{\L\ba i}^{\bar\s}[\eta_{\D\ba i}]\Big[\int_S e^{-H_i(\s_i\bar{\s}_{\L^c}\cdot_{\L\ba i})}\a(d\s_i) \Big]}
\end{equation}

On the r.h.s. we see a RCFLM $\mu_{\L\ba i}^{\bar\s}[\eta_{\D\ba i}]$ appearing 
with constrained measure $\a_{\eta_i}$ only in the volume $\D\ba i$, i.e.
\begin{eqnarray}
\mu_{\L\ba i}^{\bar\s}[\eta_{\D\ba i}](d\s_\L)=\frac{e^{-H_{\L\ba i}(\s_{\L\ba i}\bar\s_{\L^c})}\prod_{i\in\D\ba i}K(d\s_i|\eta_i)\prod_{i\in\L\ba\D}\a(d\s_i)}{\int_{S^{\L\ba i}}e^{-H_{\L\ba i}(\td\s_{\L\ba i}\bar\s_{\L^c})}\prod_{i\in\D\ba i}K(d\td\s_i|\eta_i)\prod_{i\in\L\ba\D}\a(d\td\s_i)}
\end{eqnarray} 

By the assumption of Theorem \ref{main1} we can construct the measures on the first 
layer as an infinite-volume limit with boundary condition $\bar \s$. 

Hence, the conditional  distribution 
$\g'_{i,\L,\bar\s}(d\eta_{i}|\eta_{\D\ba i})$ has an infinite-volume limit $\g'_{i,
\bar\s}(d\eta_{i}|\eta_{\D\ba i})$, for any arbitrary 
conditioning $\eta_{\D\ba i}$, since $h(\s_{\L\ba i}\bar\s_{\L^c},\eta_i):=\int_S e^{-H_i(\s_\L\bar{\s}_{\L^c})}k(\s_i,\eta_i)\a(d\s_i)$ 
is a bounded quasilocal function in $\s$ for each $\eta_i$. 
Note that this conditional distribution still depends on the boundary condition $\bar \s$ 
when the initial specification is in the phase transition regime. 
Let us denote the corresponding specification of the RCFLM with $\eta$-conditioning only in $\D\ba i$  
 by $\g[\eta_{\D\ba i}]$. It follows from (\ref{last}) that
\begin{equation}\begin{split} 
&\Big\Vert\g'_{i,\bar\s}(d\eta_{i}|\eta_{\D\ba i})-\g'_i(d\eta_{i}|\eta_{i^c})\Big\Vert\cr
&\leq
2\frac{\Big|\mu_{i^c}[\eta_{\D\ba i}]\Big[\int_{S'}h(\cdot,\eta_i)\a'(d\eta_i)\Big]-\mu_{i^c}[\eta_{i^c}]\Big[\int_{S'}h(\cdot,
\eta_i)\a'(d\eta_i)\Big]\Big|}
{\mu_{i^c}[\eta_{i^c}]\Big[\int_{S'}h(\cdot,\eta_i)\a'(d\eta_i)\Big]}.
\end{split}
\end{equation}

But  using   the fact that
\begin{equation}
\Vert \g_j[\eta_{\D\ba i}] -\g_j[\eta_{i^c}] \Vert  \left\{ \begin{array}{ll}
         =0 & \mbox{if $j\in\D\ba i $};\\
        \leq2 & \mbox{if $j\in \D^c$}\end{array} \right.  
\end{equation}
we have by the comparison criterion in Proposition \ref{Dobdist} and  using the  assumption that the  RCFLM
with full $\eta$-conditioning  satisfies  Dobrushin's condition  uniformly in $\eta$ that
\begin{equation}\begin{split}
&\Big|\mu_{i^c}[\eta_{\D\ba  i}]\Big[\int_{S'}h(\cdot,\eta_i)\a'(d\eta_i)\Big]-\mu_{i^c}[\eta_{i^c}]\Big[\int_{S'}h(\cdot,\eta_i)\a'(d\eta_i)
\Big]\Big|\cr 
& \leq2\sum_{i\in G}\sum_{j\in \D^c}\d_i\Big(\int_{S'}h(\cdot,\eta_i)\a'(d\eta_i)\Big)\bar D_{ij}.
\end{split}
\end{equation}
Taking now the limit $\D\uparrow G$ we get (\ref{gamap}), by weak  convergence of the RCFLM in $\D$
to the full one, and by the backwards martingale convergence theorem. The convergence is weak since we require
the single spin space to be separable and metrizable. In this set-up weak quasilocal topology is equivalent 
to weak topology.

\item
The proof of assertion 2 utilizes  the definition of the single-site
part of the RCFLM  and arbitrary test function  $g$, with $|g|\leq
1$ to   define
\begin{equation}\label{ker}
\Bigg|\int g(\s_j)\Bigl(
 \g_j[\eta_{j}](d\s_j|\bar\s_{G\ba
i})-\mu_j[\bar\eta_{j}]
 (d\s_j|\bar\s_{G\ba i})\Bigr)\Bigg|.
\end{equation}
The rest of the proof follows  by adding and subtracting the
following quantity
\begin{equation}\label{int}
\dfrac{\int g(\s_j)\exp\Bigl(-H_j(\s_j\bar\s_{G\ba\{i,j\}})
\Bigr)K(d\s_j|\eta_j)\int
\exp\Bigl(-H_j(\td\s_j\bar\s_{G\ba\{i,j\}})
\Bigr)K(d\td\s_j|\eta_j)}{\int
\exp\Bigl(-H_j(\td\s_j\bar\s_{G\ba\{i,j\}})
\Bigr)K(d\td\s_j|\bar\eta_j)\int
\exp\Bigl(-H_j(\td\s_j\bar\s_{G\ba\{i,j\}})
\Bigr)K(d\td\s_j|\eta_j)}
\end{equation}
 to the expression under the absolute value sign in
 (\ref{ker}), rearranging terms and simplifying appropriately.

 \item It follows from (\ref{38}) and (\ref{39}) of Proposition
 \ref{Dobdist} that
 \begin{equation}\begin{split}
&\Big|\mu_{i^c}[\eta_{i^c}](h_2)-\mu_{i^c}[\bar\eta_{i^c}](h_2)\Big|\cr
& \leq 2\sum_{k,j\in G\ba i}\d_k(h_2)\bar D_{kj}
\exp\Bigl(\sum_{A\ni j}\d_j(\P_A)\Bigr)\Big\Vert
K(\cdot|\eta_j)-K(\cdot|\bar\eta_j)\Big\Vert,
\end{split}
\end{equation}
 since by definition of $H_i$,  $h_2$ is  a local
 function on $S^{G\setminus i}$. The rest of the proof of 3 follows
 from the  bound on  $\d_k(h_2)$  given in statement 4 of the
 Proposition.

 \item Recalling that $h_2(\s_{i^c})=\int_{S\times
S'}K(d\s_{i},d\eta_{i})\exp\Bigl(-H_{i}(\s_{i}\s_{i^c})\Bigr)$ we
estimate for any pair of configurations $\s$ and $\bar \s$ that
coincide except on $k$
\begin{equation}\begin{split}
&\Bigl|
\exp\Bigl(-H_{i}(\s_{i}\s_{i^c})\Bigr)-\exp\Bigl(-H_{i}(\s_{i}\bar\s_{i^c})\Bigr)
\Bigr| \cr &= \Bigl| \exp\Bigl(-\sum_{A\ni
i,k}\P_A(\s_{i}\s_{i^c})\Bigr)- \exp\Bigl(-\sum_{A\ni
i,k}\P_A(\s_{i}\bar\s_{i^c})\Bigr) \Bigr| \exp\Bigl(-\sum_{A\ni i,
A\not\ni k}\P_A(\s_{i}\s_{i^c})\Bigr)\cr &\leq \d_k\Bigl( \sum_{A\ni
i,k}\P_A\Bigr) e^{   \sum_{A\ni i}\Vert \P_A\Vert_\infty},
\end{split}
\end{equation}
where we have used the fact that $|e^x-e^y|\leq
|x-y|e^{\max\{x,y\}}$.
 \item Take a test function $\f:S'\rightarrow \R$, with
$|\f|\leq 1$ and  consider
\begin{equation}\begin{split} \label{gdist}
&\int_{S'}\f(\eta_{i})\Bigl(\g'_i(d\eta_{i}|\eta_{i^c})-
\g'_i(d\eta_{i}|\bar\eta_{i^c})\Bigr)\cr &=\frac{\int_{S^{G\ba
i}}\mu_{i^c}[\eta_{i^c}](d\s_{i^c})h_1(\s_{i^c})}{\int_{S^{G\ba
i}}\mu_{i^c}[\eta_{i^c}](d\s_{i^c})h_2(\s_{i^c})}
 -\frac{\int_{S^{G\ba
 i}}\mu_{i^c}[\bar\eta_{i^c}](d\s_{i^c})h_1(\s_{i^c})}{\int_{S^{G\ba
 i}}\mu_{i^c}[\bar\eta_{i^c}](d\s_{i^c})h_2(\s_{i^c})},
\end{split}
\end{equation}
where we have set  $h_1(\s_{i^c})=\int_{S\times
S'}K(d\s_{i},d\eta_{i})\f(\eta_{i})\exp\Bigl(-H_{i}(\s_{i}\s_{i^c})\Bigr)$
. By adding and subtracting\\
$\dfrac{\mu_{i^c}[\eta_{i^c}](h_1)\mu_{i^c}[\eta_{i^c}](h_2)}{\mu_{i^c}
[\eta_{i^c}](h_2)\mu_{i^c}[\bar\eta_{i^c}](h_2)}$ to the right hand
side of (\ref{gdist}) and making use of the fact that
$||\f||_\infty\leq 1$ yields

\begin{equation}\begin{split}\label{gd}
&\Big|\int_{S'}\f(\eta_{i})\Bigl(\g'_i(d\eta_{i}|\eta_{i^c})-\g'_i
(d\eta_{i}|\bar\eta_{i^c})\Bigr)\Big| \leq
2\frac{\Big|\mu_{i^c}[\eta_{i^c}](h_2)-\mu_{i^c}[\bar\eta_{i^c}]
(h_2)\Big|}{\mu_{i^c}[\bar\eta_{i^c}](h_2)}.
\end{split}
\end{equation}
 \end{enumerate}
\begin{flushright}
$\Cox$
\end{flushright}
Note from the  proof of statement 5 of the above Proposition   that
the denominator in (\ref{gd}) can as well be
$\mu_{i^c}[\eta_{i^c}](h_2)$ if one adds and subtracts from the
right hand side of (\ref{gdist})
$\dfrac{\mu_{i^c}[\bar\eta_{i^c}](h_1)\mu_{i^c}[\bar\eta_{i^c}](h_2)}
{\mu_{i^c}[\eta_{i^c}](h_2)\mu_{i^c}[\bar\eta_{i^c}](h_2)}$ instead
of
$\dfrac{\mu_{i^c}[\eta_{i^c}](h_1)\mu_{i^c}[\eta_{i^c}](h_2)}{\mu_{i^c}
[\eta_{i^c}](h_2)\mu_{i^c}[\bar\eta_{i^c}](h_2)}$, as was the case
in the above proof. But any of the two makes no difference since in
our estimate we don't make use of the actual integral of $h_2$ but
instead we utilize its uniform  norm.
Having disposed of the  results above, we now return to the\\\\
\textbf{Proof of Theorem \ref{main1}:}
\begin{enumerate}
\item
The proof follows from  Lemma \ref{kenel} and the unicity of the
Gibbs measures admitted by the RCFLM, which is  uniform in
$\eta$.
\item
Using  (\ref{last}) and (\ref{deno3}) of Proposition \ref{37} we 
get
\begin{equation}\begin{split}
&\Big\Vert\g'_i(d\eta_{i}|\eta_{i^c})-\g'_i(d\eta_{i}|\bar\eta_{i^c})\Big\Vert
\leq
2\frac{\Big|\mu_{i^c}[\eta_{i^c}](h_2)-\mu_{i^c}[\bar\eta_{i^c}](h_2)
\Big|}{\mu_{i^c}[\bar\eta_{i^c}](h_2)} \cr &\leq 4 e^{2\sum_{A\ni
i}\Vert \P_A\Vert_\infty} \sum_{k,j\in G\ba
i}\d_k\Big(\sum_{A\supset \{i,k\}}\P_A\Big) \Bar D_{kj}
e^{\sum_{A\ni j}\d_j(\P_A)}\Big\Vert
K(\cdot|\eta_j)-K(\cdot|\bar\eta_j)\Big\Vert.
\end{split}
\end{equation}
The 2 in front of $\sum_{A\ni i}\Vert \P_A\Vert_\infty$ in the
 exponential  is obtained by observing that
 $\frac{1}{\mu_{i^c}[\bar\eta_{i^c}](h_2)}\leq \frac{1}{e^{-\sum_{A\ni
 i}\Vert \P_A\Vert_\infty}}$.
 \end{enumerate}
$\Cox$

\section{ Proof of results on  short-time Gibbsianness for
 time-evolved rotator models}
\textbf{Proof of Theorem \ref{thm:main-intro}:} Consider the rotator
model on the  lattice   $G$, with $S=S^{q-1}$ ( the sphere in
q-dimensional Euclidean space, with $q\geq2$) as the spin  space and
Hamiltonian  given by $H(\s)=\sum_{i,j\in G;i\neq
j}J_{ij}\s_i\cdot\s_j$. We consider the RCFLM for this Hamiltonian
with   $K$  given by the diffusion or the heat kernel $k_t$ on the
sphere, i.e.
$K(d\s_i,d\eta_i)=K_t(d\s_i,d\eta_i)=k_t(\s,\eta)\a_o(d\s)\a_o(d\eta)$,
where $\a_o$ is the equidistribution on  $S^{q-1}$. In this case we
have $S=S'=S^{q-1}$  and  $\a_o(d\s_i)=\int_{\eta_i}
K_t(d\s_i,d\eta_i)$. For  the given Hamiltonian,
$H_{i}(\cdot\z_{i^c})-H_{i}(\cdot\bar\z_{i^c})$ is Lipschitz
continuous with
 Lipschitz constant
$L_{ij}=2|J_{ij}|$. To obtain the desired bound on the Dobrushin
interdependence matrix entries we employ the bound given by Corollary
\ref{maincor}. In view of this,  we need to evaluate  the integrals
$\int_S d^2(\s_i,a_i) K(d\s_i|\eta_i)=\int_{S^{q-1}}
d^2(\s_i,a_i)k_t(\s_i,\eta_i)\a_o(d\s_i) $. To compute this
integrals we choose $a_i=\eta_i$ and  denote by $Z^q_t$   the $q$-th
coordinate of a diffusion  on the sphere started at $Z^q_{t=0}=1$ (
the "north-pole") and denote the corresponding expectation by $\E$.
Thus for any $\eta_i$ we have;
\begin{equation}\begin{split}\label{ehoch}
 &  \int \a_o(d\s_i) k_{t}(\s_i,\eta_i) d^2(\s_i,\eta_i) = 2 (1- \E
 Z^q_t
 )=  2 (1- e^{-(q-1)t}).
\end{split}
\end{equation}
The first equality uses the idea that Brownian motion on the sphere
is rotation invariant and consequently choosing
$\eta_i=(0,\cdots,0,1)$.
 To see the last equality use either an explicit form of
 the transition kernel $k_t$ in polar coordinates and orthogonality
 of Legendre polynomials as in \cite{MUL}. Or use  that the generator
 of the diffusion $Z^q_t$  given by the $u$-dependent
 parts of the Laplace-Beltrami operator on the sphere  reads
 $(1-u^2)\bigl(\frac{d}{du}\bigr)^2
 -(q-1)u \frac{d}{du}$ and generates the equation
 $\frac{d}{dt}\E Z^q_t = -(q-1)\E Z^q_t $. Solving with the initial
 condition $Z^q_{t=0}=1$ yields the desired result. Note  in our
 present
 set-up that for any pair $i,j\in G\ba i_o$ we have
 $\sum_{A\supset\{i,j\},i_o\notin A }\d(\P_A)=2|J_{ij}|.$
Then it follows from Corollary \ref{maincor} that
\begin{equation}\begin{split}
&c'[\eta]\leq \sqrt2\sup_{i_o\in G}\sup_{i\in G\ba
 i_o}\sum_{j\in
G\sm i_o}\exp\bigr( |J_{ij}|\bigl)|J_{ij}|\left(1-e^{-(q-1)t}
\right)^{\frac{1}{2}} \cr &\leq \sqrt2\sup_{i\in G}\sum_{j\in
 G}e^{|J_{ij}|}|J_{ij}|\,\left(1-e^{-(q-1)t}
\right)^{\frac{1}{2}}.
\end{split}
\end{equation}
The above estimate on $c'[\eta]$ is  uniform in
 $\eta$. \\\\
1.$\;$Therefore  the proof of the  Gibbsianness of the time-evolved
measure $\mu_t$ follows from the above uniform estimate on $c'[\eta]$ 
and the hypothesis of the theorem. \\\\
2.$\;$An application of the continuity estimate on $\g'_i$ in Theorem
 \ref{main1} to the rotator model yields a continuity estimate on
 $\g'_{i,t}$ when we define $Q_{ij}(t)$ by the bound on $Q_{i,j}$ in
 (\ref{qij}). Since  the introduction
 of the Euclidean metric $d$ follows from the estimate on  the
 posterior metric $d'$
 found in Proposition \ref{Festimates} and  the quantity $\bar D_{kj}$
 appearing in the
 definition of  the $Q_{ij}$ in Theorem \ref{main1} is given by
 $\bar D_{kj}(t)=\big(\mathbf{1}+ \sum_{n=1}^\infty
 \big(1-e^{-(q-1)t}\big)^{\frac{n}{2}} A^n\big)_{kj}$ where
 $A_{ij}=e^{|J_{ij}|}|J_{ij}|$. It is also
 elementary to see that $\sum_{A\ni j}\d_j(\P_A)\leq 2\sum_{A\ni
 j}||\P_A||_\infty$  and for each $i\in G$,  $\sum_{A\ni
 i}||\P_A||_\infty=\sum_{j\in G}|J_{ij}|$. Thus, putting all the above 
 together we get
 \begin{equation}\begin{split}
&\big\|\g'_{i,t}(\cdot|\eta_{i^c})-\g'_{i,t}(\cdot|\bar\eta_{i^c})
\big\|\leq 
\sqrt{\frac{\pi}{4 t}}
\sum_{j\in G\ba
i}Q_{ij}(t)d(\eta_j,\bar\eta_j).
\end{split}
\end{equation}
The rest of the proof follows from a telescoping argument involving
the sites in $G\ba i$. The main result in this direction that we
will employ in our proof is formulated in the lemma below.

\begin{lem}\label{mainlemma}
For each non-empty finite subset $V_1\subset G\ba i$ we have the
following estimate
\begin{equation}\begin{split}
&\big\|\g'_{i,t}(\cdot|\eta_{i^c})-\g'_{i,t}(\cdot|\bar\eta_{i^c})
\big\|\leq \frac{1}{2}\sum_{j\in
V_1}\min\bigg\{
\sqrt{\frac{\pi}{t}}
Q_{ij}(t),e^{4\sum_{k\in
G}|J_{jk}|}-1\bigg\}d(\eta_j,\bar\eta_j)\cr &+
\big\|\g'_{i,t}(\cdot|\eta_{V_1^c\ba
i}\bar\eta_{V_1})-\g'_{i,t}(\cdot|\bar\eta_{ i^c}) \big\|.
\end{split}
\end{equation}
\end{lem}
Note from the second term in  the above bound that the conditionings
coincides in the chosen finite volume $V_1$. We proceed by appling the
Lemma \ref{mainlemma} to obtain a similar bound for $\big\|\g'_{i,t}(\cdot|\eta_{V_1^c\ba
i}\bar\eta_{V_1})-\g'_{i,t}(\cdot|\bar\eta_{ i^c}) \big\|$ this time for
any non-empty finite subset  $V_2\subset G\ba V_1\cup\{i\}$. Thus we
have
\begin{equation}\begin{split}
&\big\|\g'_{i,t}(\cdot|\eta_{i^c})-\g'_{i,t}(\cdot|\bar\eta_{i^c})
\big\|\leq \frac{1}{2}\sum_{j\in V_1\cup
V_2}\min\bigg\{
\sqrt{\frac{\pi}{t}}
Q_{ij}(t),e^{4\sum_{k\in
G}|J_{jk}|}-1\bigg\}d(\eta_j,\bar\eta_j)\cr &+
\big\|\g'_{i,t}(\cdot|\eta_{(V_1\cup V_2)^c\ba i}\bar\eta_{V_1\cup
V_2})-\g'_{i,t}(\cdot|\bar\eta_{ i^c}) \big\|.
\end{split}
\end{equation}
Successive application of Lemma \ref{mainlemma} along such sequence
of pair-wise disjoint non-empty finite subsets $V_n$ such that $\cup_n V_n=G\ba i$ yields the desired
result.
\begin{flushright}
$\Cox$
\end{flushright}
\textbf{Proof of Lemma \ref{mainlemma}:}\\
For any non-empty finite subset $\L\subset G\ba i$ we let
$n_\L:\L\longrightarrow \{1,2,\cdots, |\L|\}$ be a bijection between
$\L$ and $\{1,2,\cdots, |\L|\}$ and denote by $\bar\eta_{l\leq}\eta$
the configuration that coincides with $\bar\eta$ on
$n_{\L}^{-1}\big(\{1,\cdots,l\}\big)$ and $\eta$ on $G\ba
n_{\L}^{-1}\big(\{1,\cdots,l\}\big)\cup\{i\}$. The map $n_\L$ orders
the elements in $\L$. For $G=Z^2$ this map can be a spiral ordering
of the sites in $\L$. Recall that the joint a priori measure
$K_t(d\s_i,d\eta_i)
 =k_t(\s_i,\eta_i)\a_o(d\s_i)\a_o(d\eta_i)$ where as before $\a_o=\int
 K_t(\cdot,d\s_i)$.
  In this way we can write the single-site part of $\g'$ as;
 \begin{equation} \begin{split}\label{fi}
  &
 \g'_{i,t}(d\eta_i|\eta_{i^c})=f(\eta_i|\eta_{i^c})\a_o(d\eta_i),\quad
  \text{where}\cr
  &f(\eta_i|\eta_{i^c})=\dfrac{\int_{S^{G\ba
 i}}\mu_{i^c}[\eta_{i^c}](d\td{\s}_{i^c})
  \int_{S}\exp\Big(-H_i(\s_i\td{\s}_{i^c})\Big)k_t(\s_i,\eta_i)\a_o(d\s_i)}
  {\int_{S^{G\ba i}}\mu_{i^c}[\eta_{i^c}](d\td{\s}_{i^c})
  \int_{S}\exp\Big(-H_i(\s_i\td{\s}_{i^c})\Big)\a_o(d\s_i)}.
 \end{split}
 \end{equation}
  With the order on $\L$ we can now write for any pair of
  conditionings $\eta,\bar\eta\in\O'=(S')^G$
  \begin{eqnarray}\label{teles}
  f(\eta_i|\eta_{i^c})-f(\eta_i|\bar{\eta}_{i^c})=\sum_{l=1}
  ^{|\L|+1}\nabla_{l}f(\eta_i|
  \eta_{i^c},\bar{\eta}_{i^c})\quad \text{with}\\
  \end{eqnarray}
 \[ \nabla_{l}f(\eta_i|
  \eta_{i^c},\bar{\eta}_{i^c}) = \left\{ \begin{array}{ll}
  f(\eta_i|\bar\eta_{l-1\leq}\eta)-f(\eta_i|\bar\eta_{l\leq}\eta)
 & \mbox{if $1 \leq l\leq |\L|$};\\
 f(\eta_i|\bar\eta_{|\L|\leq}\eta)-f(\eta_i|\bar\eta_{i^c})
 & \mbox{if $ l=|\L|+1$,}
 \end{array} \right. \]

 where we assume $\{1,\cdots,l-1\}=\emptyset$ for $l=1$.
 In this spirit it follows from the triangle inequality that
 \begin{equation}\begin{split}
 &\big\|\g'_{i,t}(\cdot|\eta_{i^c})-\g'_{i,t}(\cdot|\bar\eta_{i^c})\big\|=\int\a_o(d\eta_i)\bigg|
 \sum_{l=1}^{|\L|+1}\nabla_l
 f(\eta_i|\eta_{i^c},\bar{\eta}_{i^c}) \bigg|\cr
 &\leq\sum_
 {l=1}^{|\L|}
 \int_{S'}\a_o(d\eta_i)\big|\nabla_l
 f(\eta_i|\eta_{i^c},\bar{\eta}_{i^c})\big|+\big\|\g'_{i,t}(\cdot|\bar\eta_{\L}\eta_{G\ba \{i\}\cup \L})-\g'_{i,t}(\cdot|\eta_{i^c})\big\|.
 \end{split}
 \end{equation}
To get the desired bound for the first term  in the above inequality we use  two
estimation procedures which provide bounds for the terms in the sum that are multiples 
of $d(\eta_j,\bar\eta_j)$. 
 
As a first step  we consider for any $1\leq l\leq |\L|$
an estimate similar to the one given in (\ref{deno}) but here we
define $Q_{ij}(t)$ by the bound in (\ref{qij}). Note for $1\leq
l\leq |\L|$ the conditionings in the definition of $\nabla_l
f(\cdot|\eta_{i^c},\bar\eta_{i^c})$ coincide except at the site
$j=n_\L^{-1}(l)$. Thus it follows from (\ref{deno}) and the estimate
on the posterior metric in Proposition \ref{Festimates} that   for
each $1\leq l\leq |\L|$
\begin{equation}\begin{split}\label{estimate1}
 &\big\|\g'_{i,t}(\cdot|\bar\eta_{l-1\leq}\eta)-\g'_{i,t}(\cdot|\bar\eta_{l\leq}\eta)\big\|=
 \int_{S'}\a_o(d\eta_i)\big|\nabla_l
 f(\eta_i|\eta_{i^c},\bar{\eta}_{i^c})\big|\leq
  \sqrt{\frac{\pi}{4 t}}
 Q_{ij}(t)d(\eta_j,\bar{\eta}_j).
\end{split}
\end{equation}

Next we apply the following estimation technique to obtain a
second bound on\\
$\big\|\g'_{i,t}(\cdot|\bar\eta_{l-1\leq}\eta)-\g'_{i,t}(\cdot|\bar\eta_{l\leq}\eta)\big\|$
for $1\leq l\leq |\L|$. First we   set $j=n_\L^{-1}(l)$  and note
that
\begin{equation}\begin{split}\label{qfi}
&\frac{f(\eta_i|\bar\eta_{l-1\leq}\eta)}{f(\eta_i|\bar\eta_{l\leq}\eta)}=\frac{f(\eta_j|\eta_i\bar\eta_
{l-1\leq}\eta_{l>})}
{f(\bar\eta_j|\eta_i\bar\eta_{l-1\leq}\eta_{l>})}\times
 \frac{\int_{S'}f(\eta_i,\bar\eta_j|\bar\eta_{l-1\leq}\eta_{l>})\a_o(d\eta_i)}
{\int_{S'}f(\eta_i,\eta_j|\bar\eta_{l-1\leq}\eta_{l>})\a_o(d\eta_i)},
\end{split}
\end{equation}
where $\bar\eta_{l-1\leq}\eta_{l>}$ is the configuration that coincides
 with $\bar\eta$ on $n_{\L}^{-1}\big(\{1,\cdots,l-1\}\big)$ and $\eta$
 on $G\ba n_{\L}^{-1}\big(\{1,\cdots,l-1\}\big)\cup\{i,j\}$ and
  $f(\eta_i,\eta_j|\bar\eta_{l-1\leq}\eta_{l>})$ is given by (\ref{fi}) if we
 appropriately replace $i$ in (\ref{fi}) with $\{i,j\}$.\\
Therefore  setting $h_2(\s_{j^c},\eta_j)=\int_{S
 }\exp(-H_j(\s_j\s_{j^c}))k_t(d\s_j,\eta_j)\a_o(d\s_i)$
we have
\begin{equation}\begin{split}\label{70}
&\frac{f(\eta_j|\eta_i\bar\eta_{l-1\leq}\eta_{l>})}{f(\bar\eta_j|\eta_i\bar\eta_{l-1
leq}\eta_{l>})}=\frac{\mu_{j^c}[\eta_i\bar\eta_{l-1\leq}\eta_{l>}]\Big
[h_2(\s_{j^c},\eta_j)\Big]}{\mu_{j^c}[\eta_i\bar\eta_{l-1\leq}
\eta_{l>}]\Big[h_2(\s_{j^c},\bar\eta_j)\Big]}.
\end{split}
\end{equation}
Let $R$ be a rotation such that $R \bar\eta_{j}=\eta_j$ and set
$\s'_j=R\s_j$. Then  it follows from the fact that
$|H_j(\s_j\s_{j^c})-H_j(\s'_j\s_{j^c})|\leq \Big(\sum_{k\in
G}|J_{jk}|\Big)d(\eta_j,\bar\eta_j)$
\begin{equation}\begin{split}
& h_2(\s_{j^c},\eta_j)=\cr &\int_S\bigg\{\int_S
\exp\bigg(-\Big(H_j(\s_j\s_{j^c})-H_j(\s'_j\s_{j^c})\Big)-H_j
(\s'_j\s_{j^c})
\bigg)K_t(d\s'_j|\eta_j)\bigg\}K_t(d\s_j|\eta_j)\cr &\leq
\exp\Big(c_j
\;d(\eta_j,\bar\eta_j)\Big)\int_{S}\exp\Big(-H_j(\s'_j\s_{j^c})
\Big)K_t(d\s'_j|\eta_j)\quad \text{and \; similarly}\cr 
& h_2(\s_{j^c},\bar\eta_j)\leq  \exp\Big(c_j
\;d(\eta_j,\bar\eta_j)\Big)\int_{S}\exp\Big(-H_j(\s_j\s_{j^c})
\Big)K_t(d\s_j|\eta_j)
\end{split}
\end{equation}
where $c_j=\sum_{k\in G}|J_{jk}|.$ It follows from (\ref{70}) and
the rotation invariance of $K_t$ that
\begin{equation}\begin{split}
&\frac{f(\eta_j|\eta_i\bar\eta_{l-1\leq}\eta_{l>})}{f(\bar\eta_j|
\eta_i\bar\eta_{l-1\leq}\eta_{l>})}\leq\frac{\mu_{j^c}[\eta_i\bar\eta_
{l-1\leq}\eta_{l>}]\Big[\exp\Big(c_j
\;d(\eta_j,\bar\eta_j)\Big)\int_{S}\exp\Big(-H_j(\s'_j\s_{j^c})
\Big)K_t(d\s'_j|\eta_j)\Big]}
 {\mu_{j^c}[\eta_i\bar\eta_{l-1\leq}\eta_{l>}]\Big[\int_{S
 }\exp(-H_j(\s'_j\s_{j^c}))k_t(d\s'_j,\eta_j)\a_o(d\s_i)\Big]}\cr
 & = e^{c_jd(\eta_j,\bar\eta_j)}.
\end{split}
\end{equation}
The above estimate follows by applying the rotation $R$ to the $\bar\eta_j$ in the r.h.s. of (\ref{70}).
Furthermore, it is not hard to deduce that
\begin{equation}\begin{split}
&\frac{\int_{S'}f(\eta_i,\bar\eta_j|\bar\eta_{l-1\leq}\eta_{l>})\a_o(d\eta_i)}
{\int_{S'}f(\eta_i,\eta_j|\bar\eta_{l-1\leq}\eta_{l>})\a_o(d\eta_i)}
\leq\sup_{\eta_i}\frac{f(\bar\eta_j|\eta_i\bar\eta_{l-1\leq}\eta_{l>})}{f(\eta_j|
\eta_i\bar\eta_{l-1\leq}\eta_{l>})}\leq e^{c_j d(\eta_j,\bar\eta_j)}.
\end{split}
\end{equation}
Therefore it follows from (\ref{qfi})  that
\begin{equation}\begin{split}
&\frac{f(\eta_i|\bar\eta_{l-1\leq}\eta)}{f(\eta_i|\bar\eta_{l\leq}\eta)}
\leq e^{2c_jd(\eta_j,\bar\eta_j)}.
\end{split}
\end{equation}
Hence  for any  $1\leq l\leq |\L|$ we have
\begin{equation}\begin{split}\label{estimate2}
&\int_{S}\a_o(d\eta_i)\bigg|\nabla_l
f(\eta_i|\eta_{i^c},\bar\eta_{i^c})\bigg|=
\int_{S}\a_o(d\eta_i)\bigg|\Big(\frac{f(\eta_i|\bar\eta_{l-1\leq}\eta)}
{f(\eta_i|\bar\eta_{l\leq}\eta)}-1\Big)f(\eta_i|\bar\eta_{l\leq}\eta)
\bigg|\cr & \leq e^{2c_jd(\eta_j,\bar\eta_j)}-1\leq
\frac{e^{4c_j}-1}{2}d(\eta_j,\bar\eta_j).
\end{split}
\end{equation}

Comparing (\ref{estimate1}) and (\ref{estimate2}) it is clearly seen
 for any $1\leq l\leq|\L|$ with $j=n_\L^{-1}(l)$ that
\begin{equation}
\big\|\g'_{i,t}(\cdot|\bar\eta_{l-1\leq}\eta)-\g'_{i,t}
(\cdot|\bar\eta_{l\leq}\eta)\big\|\leq \frac{1}{2}
\min\Big\{
\sqrt{\frac{\pi}{t}}
Q_{ij}(t),e^{4\sum_{k\in
G}|J_{jk}|}-1\Big\}d(\eta_j,\bar\eta_j),
\end{equation}
which proves the lemma.
 $\Cox$ \\

Lemma \ref{mainlemma} has an extension for interactions for which
$H_j(\cdot\s_{j^c})$ is not Lipschitz continuous. In
this set-up we have for any non-empty finite subset $V\subset G\ba
i$
 \begin{equation}\begin{split}
&\big\|\g'_{i,t}(\cdot|\eta_{i^c})-\g'_{i,t}(\cdot|\bar\eta_{i^c})
\big\|\leq \sum_{j\in V}\bigg(e^{4\d_j\Big(\sum_{A\ni
j}\P_A\Big)}-1\bigg) + \big\|\g'_{i,t}(\cdot|\eta_{V^c\ba
i}\bar\eta_{V})-\g'_{i,t}(\cdot|\bar\eta_{ i^c}) \big\|.
\end{split}
\end{equation}

To obtain the desired bound on the posterior metric we need to solve
 the diffusion equation
on the sphere $S^{q-1}$. However, it turns out in the analysis that
we
 don't need all
the components of the diffusion to arrive at our desired bound.  The
 only coordinate that
we will be interested in, is the $q$th, i.e.  we only have to solve the
 resulting diffusion
equation for the $q$th component. We employ both analytical and
 stochastic differential
equation (sde) techniques to arrive at the diffusion of interest. It
 turns out that the sde
approach easily provides the desired bound. Nevertheless, we present
 the
analytical approach because of its interest per se. We first state
the
 corresponding sde result.

\begin{lem}\label{F21}
\begin{enumerate}

\item Denote by  $Z_t$ the $q$th-component
of  the diffusion on the sphere $S^{q-1}$ for $q\geq 2$, started at a value $\sin \phi_0$ with $\phi_0\in (0,\frac{\pi}{2})$.  
Then there is a coupling of $Z_t$ to a Brownian motion on the line, $B_t$ 
such that the first passage time of $Z_t$ at zero, denoted by $T_0(Z)$  is dominated from 
above by that of $ \phi_0 + \sqrt{2}B_t$.  

\item Consequenty, independently of the dimension $q-1$ there is the estimate  
\begin{equation}\begin{split}\label{a13}
P(T_0(Z) \geq t) \leq  P(T_0(\phi_0 + \sqrt{2}B_\cdot)\geq t)\leq  2 P\Bigl(0\leq G\leq \frac{\phi_0}{\sqrt{ 2 t
}}\Bigr)
\end{split}
\end{equation}
where $G$ is a standard normal variable.

\end{enumerate}

\end{lem}

{\bf Proof: } Consider the case $q\geq 3$ first. 
The sde for the $q$-th component reads,
\begin{equation}\begin{split}
d Z_t&= -(q-1)Z_t dt +\sqrt{2(1-Z_t^2)}d B_t \cr
\end{split}
\end{equation}
Consider the transformation
\begin{equation}\begin{split}
Z_t=\sin(\bar \phi_t)
\end{split}
\end{equation}
to an unknown function $\bar \phi_t$ describing the elevation above the equator. 
We apply this transformation only for $0<Z_t<1$, and so 
there is a one-to-one map to $0<\bar \phi_t<\frac{\pi}{2}$. 
In this range the sde is equivalent to
\begin{equation}\begin{split}\label{phisde}
d \bar \phi_t&= -(q-2) \tan \bar \phi_t  + \sqrt{2}\, d B_t\cr 
\end{split}
\end{equation}
Indeed, for $q\geq 3$ the diffusion $\bar \phi_t$ does not leave the interval $(-\frac{\pi}{2},\frac{\pi}{2})$, 
meaning that that, with probability one the northpole is never reached by $Z_t$. (That this 
is true can be seen by projecting $Z_t$  along the $q$-th axis, onto the $q-1$-dimensional plane.)  

Integrating from zero to $t$ we obtain from (\ref{phisde})
\begin{equation}\begin{split}
\bar \phi_t&= -(q-2)\int_{0}^t \tan \bar \phi_s ds   + \sqrt{2}\, B_t + \phi_0\cr 
\end{split}
\end{equation}
From this equality we see that as long as $\bar \phi_s\geq 0$ for all $s\in [0,t]$ 
we have the bound $\bar \phi_t \leq   \sqrt{2}\, B_t + \phi_0$. 
This shows that the first passage time of $\bar \phi_t$ is not bigger than 
that of $\sqrt{2}\, B_t + \phi_0$. 

The proof of the inequality follows from bounding 
$P(T_0(Z) \geq t) $ from above by the first passage time of 
the Brownian motion on a line, $P(  T_0(\sqrt{2}\, B_\cdot + \phi_0)  \geq t)$. 
The latter can be computed exactly by the reflection principle applied 
to standard Brownian motion, as it is well-known. 
(We will use the reflection principle also in the proof Lemma \ref{F}, applied   
to the diffusion on the sphere.) 
This gives 
rise to the estimate on the r.h.s. \\ 

That the inequality holds also in the case $q=2$ (and is a strict inequality then) 
can be seen directly without making reference to the SDE. 
We note that the paths of a diffusion on the circle are given 
by Brownian motions on the angular variable, i.e. 
$\bar \phi_t=\sqrt{2}\, B_t + \phi_0$. Then 
 $\bar\phi_t=0$ implies 
that $Z^{(2)}_t=\sin(\bar\phi_t)=0$, but the converse is not true. 

It is interesting to realize that this construction 
provides a coupling such that  
$Z^{(q)}_t\leq   \sqrt{2}\, B_t + \phi_0$, for $q\geq 3$, 
 $Z^{(2)}_t\leq   \sqrt{2}\, B_t + \phi_0$ but {\it not} 
$Z^{(q)}_t\leq   Z^{(2)}_t$. The latter relation is guaranteed to hold only 
as long as $0\leq \sqrt{2}\, B_t + \phi_0\leq \frac{\pi}{2}$. 
$\Cox$. 
\bigskip

We now present an analytical treatment for the diffusions considered
above. This involves the study of eigenvalue problem involving the
$q$th-component of the Laplace-Beltrami operator on the sphere. In
fact the resulting eigenfunctions solve the spatial part of the
$q$th-component of the diffusion on the sphere. The transition
kernel $k_t$(defined below) for the $q$th-component of the diffusion
is determined by the solution for the above mentioned eigenvalue
problem. It is known from the literature \cite{MUL} that the Legendre
polynomials constitute a complete class of eigenfunctions, i.e. the
transition kernel $k_t$ can be written in terms of the Legendre
polynomials.

\begin{defn}\label{Rod}
The Legendre polynomial $P_n(q,\cdot)$ of degree $n$ in dimension
$q\geq 2$ is given by the Rodrigues formula
\begin{equation}\begin{split}
P_n(q,s):=\dfrac{(-1)^n
\Gamma\big(\frac{q-1}{2}\big)}{2^n\Gamma\big(n+\frac{q-1}{2}\big)}\Big(1-s^2\Big)^{\frac{3-q}{2}}\Big(\dfrac
{d}{ds}\Big)^{n}\Big(1-s^2 \Big)^{\frac{q-3}{2}+n},
\end{split}
\end{equation}
where $-1\leq s\leq 1$.
\end{defn}
These Legendre  polynomials are known (see \cite{MUL} for example)
to be orthogonal and satisfy the second order differential equations
\begin{equation}
\Big[(1-s^2)\frac{d^2}{ds^2}-(q-1)s\frac{d}{ds}+n(n+q-2)\Big]P_n(q,s)=0.
\end{equation}
The last equation indicates that  the Legendre polynomials are
eigenfunctions for the eigenvalue problem for the $q$th component of
the  Laplace-Beltrami operator on the sphere $S^{q-1}$ . This
implies that the transition kernel for the $q$th coordinate $Z_t^q$
of  the Brownian motion  on  $S^{q-1}$ can be written as
\begin{equation}\begin{split}\label{kt}
&k_t(s,u):=\frac{\Gamma\big(\frac{q}{2}\big)}{\sqrt(\pi)\Gamma\big(\frac{q-1}{2}\big)}\sum_{n=0}^{\infty}
e^{-n(n+q-2)t}N(q,n)P_n(q,s)P_n(q,u),\quad\text{where}\cr 
& N(q,n) := \left\{ \begin{array}{ll}
         \frac{(2n+q-2)\Gamma(n+q-2)}{\Gamma(n+1)\Gamma(q-1)} & \mbox{if $n \geq 1$};\\
        1 & \mbox{if $n=0$}\end{array} \right. 
\end{split} 
\end{equation}
is the dimension of spherical harmonics of degree $n$ in dimension $q$. Further we have
 set $Z_{0}^q=s$ and $Z_t^q=u$, and we have also chosen the constant
$\frac{\Gamma\big(\frac{q}{2}\big)}{\sqrt(\pi)\Gamma\big(\frac{q-1}{2}\big)}$
so that for any initial $s$ the integral of $k_t(s,u)$ with respect
to the invariant  measure $(1-u^2)^{\frac{q-3}{2}}du$ (which is the
$q$-coordinate projection of the invariant surface measure on the
sphere ) over the interval [-1,1] is equal to one. We now formulate
our result on an estimate on the posterior metric
$d'(\eta_j,\bar\eta_j)$ define in (\ref{naturalmetric}). This is
given in terms of  Legendre polynomials (introduced in Definition
\ref{Rod} above) which by our construction are  also themselves
functions of $d(\eta_j,\eta'_j)$( the Euclidean distance between $\eta_j $
and $\eta'_j$ ).

\begin{lem}\label{F}
For the diffusion on a sphere there is an estimate of the
posterior-metric $d'(\eta,\eta')$ at fixed $t$ in terms of
$d(\eta,\eta')$, the induced metric on the sphere $S^{q-1}$ obtained
by imbedding the sphere into the Euclidean space, given by
\begin{equation}\begin{split}
&d'(\eta,\eta')\leq F_t( d(\eta_j,\bar\eta_j) )
 \end{split}
 \end{equation}
with the function
\begin{equation}\begin{split}
F_t( x)&= 2\Bigl(1-2 P^{\frac{x}{2}}(Z_t^q\leq 0)\Bigr)\cr
&=\dfrac{-4\Gamma\left( \frac{q}{2}\right) }{\sqrt{\pi}\Gamma\left(
\frac{q-1}{2}\right)}\sum_{n=1,3,5,\dots} e^{-n
(n+q-2)t}N(q,n)P_n\Big(q,\frac{x}{2}\Big)\int_{-1}^0
 P_n(q,s)(1-s^2)^{\frac{q-3}{2}}ds
\end{split}
\end{equation}
\end{lem}
\textbf{Proof:}\\The idea of the proof is to construct a coupling of
two diffusions on the sphere starting at the  points $\eta$ and
$\eta'$ . By rotation invariance of such diffusions we assume that
$\eta$ and $\eta'$ are mirror images of each other under reflection
at the equatorial plane. Then we construct a coupling by reflection
\cite{LIN} of the path started at $\eta$ with the equator as the
mirror line, up to the time where the diffusion hits the equator.
After that the two diffusions move on together. In this way the
coupling time  for the two diffusions is the same as  the first time
$Z_t^q=0$ (the first passage time $T_0$ to level 0 given by
$T_0:=\inf\{t\geq0, Z_t=0\}$ ) for either $Z^q_0=z$ or $Z^q_0=-z$
where $z=\varepsilon_q\cdot\eta$  (here
$\varepsilon_1,\cdots,\varepsilon_q $ constitute the canonical
orthonormal basis for $\mathbb{R}^q$  and $"\cdot"$ is the usual
scalar product ). We know from coupling theory that
\begin{equation*}
d'(\eta,\eta')\leq 2\mathbb{P}^{\frac{x}{2}}(T_0\geq t),
\end{equation*}
where $x=d(\eta,\eta')$ is the Euclidean distance between $\eta$ and
$\eta'$. Further it follows from the  reflection principle of
D\'esir\'e Andr\'e (\cite{KAR},pp.79-81 and \cite{LEV},p.293 )that
\begin{equation*}
\mathbb{P}^{\frac{x}{2}}(T_0\leq
t)=2\mathbb{P}^{\frac{x}{2}}(Z_t^q\leq 0)=2\int_{-1}^0
k_t\Big(\frac{x}{2},s\Big)(1-s^2)^{\frac{q-3}{2}}ds.
\end{equation*}
The heuristic argument for the first equality in the above equation
is as follows; the probability that the first passage time $T_0$ (to
a level $0$ for a 1-dimensional diffusion starting at some initial
point $y>0$) is less or equal to $t$ is the sum of the probabilities
of the events that  $T_0\leq t$ and $Z_t^q<0$, and  $T_0\leq t$ and
$Z_t^q>0$. The probability for the first event is the same as the
probability for the event that the 1-dimensional diffusion $Z_t^q$
starting at $y$ is below the level $0$. For the probability of the
second event observe that after the diffusion reached level $0$, it
has equal probability to reach level $-c$ below $0$ or level $c$
above $0$ since the diffusion in our set-up is symmetric about 0.
Hence the probability of the second event is the same as the first
due to the symmetry of $Z_t^q$  about 0.

It follows from  the orthogonality property of the Legendre
polynomials that for each positive even integer $n$   the integral\\
$\int_{-1}^0
P_n(q,s)(1-s^2)^{\frac{q-3}{2}}ds=\frac{1}{2}\int_{-1}^1
P_n(q,s)P_0(q,s)(1-s^2)^{\frac{q-3}{2}}ds=0$ ( since $P_0(q,s)=1$)
for all $q\geq2$. Therefore the rest of the proof follows from
(\ref{kt}) and the fact that the integral $\int_{-1}^1
P_0(q,s)^2(1-s^2)^{\frac{q-3}{2}}ds=\frac{\sqrt\pi
\Gamma\big(\frac{q-1}{2}\big)}{\Gamma\big(\frac{q}{2}\big)}.$
\begin{flushright}
$\Cox$
\end{flushright}

We have seen from the above proof that  for positive even integers
$n$ the integral  (over [-1,0] and w.r.t to the invariant  measure
$(1-s^2)^{\frac{q-3}{2}}ds$ )  of the Legendre polynomial of degree
$n$  is always equal to  zero, as long as  the dimension $q\geq2$. The
integral for the corresponding odd degree case can also be computed
explicitly and this explicit value of the integral  we formulate as
our next lemma.
\begin{lem} \label{oddL}
For any odd integer $2m+1$ (m=0,1,2,....) the integral of the
Legendre polynomials $P_{2m+1}(q,\cdot)$   over the interval [-1,0]
is given by
\begin{equation}\begin{split}
& \int_{-1}^0 P_{2m+1}(q,s)(1-s^2)^{\frac{q-3}{2}}ds=(-1)^{m}
\prod_{i=0}^m \Bigg( \dfrac{2i-1}{ q+2i-1}\Bigg).
\end{split}
\end{equation}
\end{lem}
\textbf{Proof:} We obtain from definition of $P_{2m+1}(q,s)$ in
Definition \ref{Rod} that the integral

\begin{equation}\begin{split}
& \int_{-1}^0 P_{2m+1}(q,s)(1-s^2)^{\frac{q-3}{2}}ds\cr
&=\dfrac{-1}{2^{2m+1}\prod_{i=0}^{2m}\big(2m+\frac{q-1}{2}-i\big)}
\Big(\frac{d}{ds}\Big)^{2m}\Big(1-s^2\Big)^{2m+\frac{q-1}{2}}\Big|_{s=-1}^0.\cr
\end{split}
\end{equation}
Note that  for each $m$ the  above differentiation(s) will always
involve terms which are multiples of $(1-s^2)$. This implies that
evaluating the above expression at $s=-1$ will always yield zero.
However, it follows from Binomial expansion of $\big(1-s^2\big)^r$
(where $r=2m+\frac{q-1}{2}$) that
\begin{equation}\begin{split}
&\dfrac{-1}{2^{2m+1}\prod_{i=0}^{2m}\big(2m+\frac{q-1}{2}-i\big)}\Big
(\frac{d}{ds}\Big)^{2m}\Big(1-s^2\Big)^{2m+\frac{q-1}{2}}\Big|_{s=0}\cr
&=\dfrac{(-1)^{m+1}(2m)!r (r-1)\cdots (r-(m-1))}{m!
2^{2m+1}\prod_{i=0}^{2m}\big(2m+\frac{q-1}{2}-i\big)} .
\end{split}
\end{equation}
The rest of the proof follows from the observations that $(2m)!=2^m
m!\prod_{i=1}^m (2i-1)$ and
$\frac{r(r-1)\cdots(r-(m-1))}{\prod_{i=0}^{2m}\big(2m+\frac{q-1}{2}-i\big)}
=\frac{2^{m+1}}{\prod_{i=0}^m (2i+q-1)}$.
\begin{flushright}
$\Cox$
\end{flushright}

{\bf Proof of the Proposition \ref{Festimates} :}
\begin{enumerate}
\item It follows from
 Lemma \ref{F21} that, for any $q\geq 2$,  
\begin{equation}\begin{split}
F_{q,t}(x)\leq 2\mathbb{P}^{\frac{x}{2}}\big(T_0\geq t\big)
\leq  4 P\Bigl(0\leq
G\leq  \frac{\arcsin\frac{x}{2}}{\sqrt{ 2 t }}\Bigr).
\end{split}
\end{equation}

Using $P\Bigl(0\leq G\leq u\Bigr)\leq \frac{u}{\sqrt{2 \pi}}$ by
concavity and $\arcsin y \leq \frac{\pi}{2}y$ for $0\leq y\leq 1$ we
obtain $F_{q,t}(x)\leq \frac{\sqrt{\pi}x}{2\sqrt{ t}}$. Note that in
both of the last estimates the constants were sharp.
\item
The claim for general dimensions $q\geq 2$ follows from Lemma
\ref{F} and \ref{oddL}.
 \end{enumerate}
$\Cox$

\textbf{Proof of Theorem \ref{fuzzy}:} This Theorem is an application of Theorem \ref{main1}. The only
quantities we have to worry about are the entries  of the Dobrushin interdependence matrix $\bar C$.
It follows from the hypothesis of the Theorem; namely the continuity property of the interaction  and the terms in bound on $c'[\eta]$ in Corollary \ref{maincor} that
\begin{equation}\begin{split}
 &\frac{1}{2}\exp\Bigr(
\frac{1}{2}\sum_{A\supset\{i,j\}}\delta(\Phi_A)\Bigl)L_{ij}\inf_{a_i\in S_{s'}}\Bigl(\int_{S_{s'}}
d^2(\s_i,a_i)\a_{s'}(d\s_i)\Bigr)^{\frac{1}{2}}\cr 
&\leq \sup_{s'\in S'}\frac{\rho_{s'}}{2}\exp\Bigr(
\frac{1}{2}\sum_{A\supset\{i,j\}}\delta(\Phi_A)\Bigl)L_{ij}=\bar C_{ij},
\end{split}
\end{equation}
where $\rho_{s'}:=\text{diam}(S_{s'})$ is the diameter of $S_{s'}$.
\begin{flushright}
$\Cox$
\end{flushright}
\textbf{Acknowledgements: }\\
The authors thank Aernout van Enter and Roberto Fern\'andez for interesting discussions.

\end{document}